\documentclass[a4paper,11pt,times,authoryears]{article}
\usepackage{authblk}
\usepackage[authoryear]{natbib}
\usepackage{bm}
\usepackage{amsmath}
\usepackage{amsthm}
\usepackage{mathtools}
\usepackage{amssymb,amsfonts}
\usepackage{graphicx}
\usepackage{mathrsfs}
\usepackage{multirow,multicol}
\usepackage{subcaption}
\usepackage{epstopdf}
\usepackage{enumerate}
\usepackage{booktabs}
\usepackage{adjustbox}
\usepackage{pdflscape}
\usepackage{multirow}
\usepackage[colorlinks,citecolor=blue]{hyperref}
\usepackage{rotating}
\usepackage{fullpage}
\usepackage{setspace}
\usepackage{lscape}

\usepackage{booktabs}
\usepackage{setspace}
\usepackage{tocbibind}
\usepackage{color}
\usepackage{commath} 
\usepackage{geometry}
\usepackage{setspace}
\theoremstyle{remark}
\everydisplay{\abovedisplayshortskip\belowdisplayshortskip}%
\makeatletter \oddsidemargin  -.1in \evensidemargin -.1in

\numberwithin{equation}{section}

\def\bb{\boldsymbol{\beta}}

\def\bx{{\mathbf x}}

\def\bhat{\widehat{\beta}}
\def\ahat{\widehat{\alpha}}

\def\btil{\widetilde{\beta}}
\def\atil{\widetilde{\alpha}}

\def\bI{\mathbf  I}

\def\bX{\mathbf  X}
\def\bV{\mathbf  V}
\def\bW{\mathbf  W}
\def\bZ{\mathbf  Z}

\def\mC{\mathcal{C}}
\def\mE{\mathcal{E}}
\def\mL{\mathcal{L}}

\def\aa{\alpha}
\def\bb{\beta}

\def\xi{x_{(i)}}
\def\no{\nonumber}
\def\ll{\left(}
\def\rr{\right)}

\title{\bf Estimation in Weibull Distribution Under Progressively Type-I Hybrid Censored Data}

\author[1]{Yasin Asar\thanks{Corresponding Author}}
\author[2]{R.Arabi Belaghi}
\affil[1]{Department of Mathematics--Computer, Necmettin Erbakan University, Konya, Turkey, yasar@erbakan.edu.tr, yasinasar@hotmail.com}
\affil[2]{Department of Statistics, Faculty of Mathematical Sciences, University of Tabriz, Tabriz, Iran, rezaarabi11@gmail.com}

\begin{document}
\maketitle
\bibliographystyle{agsm}
\setcounter{page}{1}
\setcounter{equation}{0}
\setcounter{figure}{0}
\setcounter{table}{0}
\noindent{\bf{Abstract:}} In this article, we consider the estimation of unknown parameters of Weibull distribution when the lifetime data are observed in the presence of progressively type-I hybrid censoring scheme. The Newton-Raphson algorithm, Expectation-Maximization (EM) algorithm and Stochastic EM (SEM) algorithm are utilized to derive the maximum likelihood estimates (MLEs) for the unknown parameters. Moreover, Bayesian estimators using Tierney-Kadane Method and Markov Chain Monte Carlo (MCMC) method are obtained under three different loss functions, namely, squared error loss (SEL), linear-exponential (LINEX) and generalized entropy loss (GEL) functions. Also, the shrinkage pre-test estimators are derived. An extensive Monte Carlo simulation experiment is conducted under different schemes so that the performances of the listed estimators are compared using mean squared error, confidence interval length and coverage probabilities. Asymptotic normality and MCMC samples are used to obtain the confidence intervals and highest posterior density (HPD) intervals respectively. Further, a real data example is presented to illustrate the methods. Finally, some conclusive remarks are presented.
\vspace{0.2cm}

\vskip 6mm
\noindent {\bf AMS SUBJECT CLASSIFICATION:} 62F10; 62N01; 62N05.

\section{Introduction}
\label{Intro}
Censored data occurs commonly in reliability and survival analysis. There are mainly two censoring schemes which are Type-I censoring where the life-testing experiment stops at a predetermined time, say $T$ and Type-II censoring, where the life-testing experiment stops when predetermined number of failures, say $m$, are observed.  \citet{epstein1954} proposed the hybrid censoring scheme which is the mixture of Type-I and Type-II censoring schemes. The hybrid censoring scheme becomes quite popular in the reliability and life-testing experiments so far. For example, see the papers of \citet{ChenBhat1988}, \citet{childs2008}, \citet{kundu2006}, \citet{BK2013}. In these schemes, it is allowed to remove the units only at the terminal points of the experiments. However, \citet{kundu2006} introduced another scheme which is called the Type-I progressively hybrid censoring scheme (Type-I PHCS) such that it allows removals of units during the test time. For more information on progressive censoring, we refer to \citet{Bala2000} and \citet{Bala2007prog}.  Type-I PHCS can be viewed as a mixture of Type-I progressive censoring and hybrid censoring as follows: Assume that there are $n$ identical units in a lifetime experiment with the progressive censoring scheme $\left(R_1, R_2, ..., R_m \right)$, $1\leq m \leq n$ and the lifetime experiment ends at a predetermined time $T \in (0, \infty)$ and $n,m, R_i$'s are all fixed integers. At the time of first failure, say $X_{1:m:n}$, $R_1$ units randomly removed from the remaining $n-1$ units. Similarly, when the second failure occurs at the time $X_{2:m:n}$, $R_2$ units are removed from the remaining $n-R_1-2$ units. This process continues up to the end of experiment which occurs at the time ${\rm min}\ll X_{m:m:n}, T\rr$. Therefore, if the $m$th failure occurs before time $T$, the experiment ends at the time $X_{m:m:n}$ and all the remaining units $R_m = n-\sum_{i=1}^{m-1}R_i -m$ are removed. However, if the experiment ends at time $T$ with only $J$ failures, $0 \leq J < m$, then all the remaining units $R_J^* = n- \sum_{i=1}^{J}R_i -J$ are removed and the test ends at time $T$. Therefore, under Type-I PHCS we have the following two cases:

\begin{itemize}
    \item 
    Case I: $\{X_{1:m:n}, X_{2:m:n}, ..., X_{m:m:n}\}$ if $X_{m:m:n} \leq T$.
    \item
    Case II: $\{X_{1:m:n}, X_{2:m:n}, ..., X_{J:m:n}\}$ if $X_{J:m:n} < T < X_{J+1:m:n}$.
\end{itemize}

Due to the fact that the lifetime distributions of many experimental units can be modeled by a two-parameter Weibull distribution which is one of the most commonly used model in reliability and lifetime data analysis, we consider the Weibull distribution in this paper. The probability distribution function (PDF) and cumulative distribution function (CDF) of two parameter Weibull distribution are given as follows:
\begin{eqnarray}\label{pdf:Weibull}
f\left(x; \alpha, \beta \right) = \alpha \beta x^{\alpha-1} {\rm exp}\left\{ -\beta x ^\alpha \right\}
\end{eqnarray}
\begin{eqnarray}\label{cdf:Weibull}
F\left(x; \alpha, \beta \right)=1-{\rm exp}\left\{-\beta x^\alpha \right\}
\end{eqnarray}
where $\alpha>0$ is the shape parameter and $\beta>0$ is the scale parameter.

\citet{BK2008} considered the statistical inference on Weibull parameters when the data are Type-II hybrid censored, maximum likelihood estimation (MLE), approximate MLE and Bayes estimation techniques were studied by the authors. \citet{Kateri2008} proposed an alternative approach based on a graphical method, which also shows the
existence and uniqueness of the MLEs. \citet{lin2009} studied the MLEs and
the approximate MLEs (AMLEs) of the parameters of Weibull distribution under adaptive Type-II progressive hybrid censoring.     \citet{HuangWu2012} discussed the maximum likelihood estimation and Bayesian estimation of Weibull parameters under progressively type-II censoring scheme. \citet{Lin2012} investigated the maximum likelihood estimation and Bayesian estimation for a two-parameter Weibull distribution based on adaptive Type-I progressively hybrid censored data which was introduced by \citet{LinHuang2012}. \citet{jia} studied
the exact inference on Weibull parameters under multiple Type-I
censoring. \citet{moktari} discussed approximate and Bayesian inferential procedures for the progressively Type-II hybrid censored data from the Weibull distribution. However, this type of censoring is identical to what we called as Type I progressive hybrid censored data.   This paper will be different from \citet{moktari} in three directions. Firstly, we introduce a new approach for inference about the Weibull distribution based on EM and SEM methods. We will show that both EM and SEM will result to have better estimates in the sense of having smaller biases and mean square errors. Secondly, We will drive the shrinkage estimators based on the ML estimates resulting to have higher deficiencies. Finally, in the Bayesian approach, different loss functions such as SEL, LINEX, and GEL will be applied with both informative and non-informative priors. 

The rest of the paper is organized as follows: In Section \ref{MLE}, MLE of the parameters are introduced by using Newton--Raphson (NR) algorithm, expectation--maximization (EM) algorithm and stochastic expectation--maximization (SEM) algorithm, also the Fisher information matrix is obtained. In Section \ref{Bayes}, Bayes estimation for the parameters of Weibull distribution under the assumption of independent priors using different loss functions such as squared error loss (SEL) function,  linear-exponential (LINEX) loss function and  general entropy loss (GEL) function. Moreover, \citeauthor{Tierney} (T-K) approximations under these loss, functions are also computed and Markov-Chain Monte Carlo (MCMC) method is also presented to estimate the parameters. In Section \ref{Shrinkage}, a shrinkage pre-test estimation method is discussed. Extensive Monte Carlo simulations are conducted and results are discussed in Section \ref{sim}. A real data example is presented in Section \ref{RealData} to illustrate the findings of the study. Finally, some conclusive remarks are given in Section \ref{conclusion}.

\section{Maximum Likelihood Estimation}\label{MLE}

Let $\bX=\left(X_{1:m:n},\ldots, X_{r:m:n}\right)$ represents the Type-I progressively hybrid censored sample of size $r$ from a sample of size $n$ drawn from a population with probability distribution given in Equation \eqref{pdf:Weibull}. Throughout this paper, we will denote $X_{i:m:n}$ by $X_{(i)}, i=1,2,\ldots,r$. Then the likelihood function of $(\alpha, \beta)$ given the observed data $\bx$ can be written as
\begin{eqnarray}
\label{likelihood:data}
L(\alpha, \beta \mid \bx) &\propto& \prod_{i=1}^{r}f(x_{(i)}; \aa, \beta)\big[1-F(x_{(i)}; \aa, \beta)\big]^{R_{i}} \big[1-F(\mathcal{C}; \aa, \beta)\big]^{R_{T}},
\end{eqnarray}
where $r=m, \mC=x_{(m)}, R_{T} =0$ in Case I, and $r=d, \mC=T, R_{T}=n-d-\sum_{i=1}^{d}R_{i}$ in Case II.  Based on the observed data, the log-likelihood function can be expressed as
\begin{eqnarray}
\label{log-likelihood:data}
l\left(\alpha, \beta \mid \bx \right)=\ln L\left(\alpha, \beta \mid \bx \right) = r\ln(\alpha \beta) +\left(\alpha-1\right)\sum_{i=1}^{r}\ln\left( x_{(i)}\right)-\beta\sum_{i=1}^{r} \left\{x_{(i)}^{\aa} \left(1+R_i\right)\right\}-\beta \mC^\alpha R_{T}.
\end{eqnarray}
Taking the derivatives of Equation \eqref{log-likelihood:data} with respect to $\alpha$ and $\beta$ and equating them to zero, one can obtain the following likelihood equations for $\alpha$ and $\beta$ respectively
\begin{eqnarray}\label{likelihood:alpha}
\frac{\partial l\left(\alpha, \beta \mid \bx \right)}{\partial \alpha} &=& \frac{r}{\alpha} + \sum_{i=1}^r \ln\left( x_{(i)}\right)-\beta \sum_{i=1}^r \left\{\left(1+R_i\right) x_{(i)}^\aa \ln\left( x_{(i)}\right)\right\} - \beta \mC^\alpha \ln(\mC) R_{T}=0\\
\label{likelihood:beta}
\frac{\partial l\left(\alpha, \beta \mid \bx \right)}{\partial \beta} &=&\frac{r}{\beta} - \sum_{i=1}^r \left\{x_{(i)}^\alpha \left(1+R_i\right)\right\} - \mC^\alpha R_{T}=0.
\end{eqnarray}
Solving Equation \eqref{likelihood:beta} yields the MLE of $\beta$ which is given by
\begin{eqnarray}\label{mle:beta}
\bhat=\frac{r}{\mC^{\ahat} R_T + \sum_{i=1}^{r} \left\{x_{(i)}^{\ahat} \left(1+R_i\right)\right\}}.
\end{eqnarray}
Now, substituting Equation \eqref{mle:beta} into \eqref{likelihood:alpha}, the MLE of $\alpha$ can be obtained by solving the following nonlinear equation:
\begin{eqnarray*}\label{mle:alpha:nonlin}
\frac{r}{\ahat}+\frac{r\left[\sum_{i=1}^{r} \left\{\left(1+R_i\right)x_{(i)}^{\ahat}\ln(\xi) \right\}+R_T\mC^{\ahat} \ln(\mC)\right]}{R_T\mC^{\ahat} +\sum_{i=1}^{r} \left\{x_{(i)}^{\ahat} \left(1+R_i\right)\right\}}=0.
\end{eqnarray*}

The second partial derivatives of the log-likelihood equation are obtained as follows:
\begin{eqnarray}\label{second:deri:alpha}
\frac{\partial^2 l\left(\alpha, \beta \mid \bx \right)}{\partial \alpha^2} &=& -\frac{r}{\alpha^2} - \beta\sum_{i=1}^r \left\{ \left(1+R_i\right) x_{(i)}^\alpha \ln\left( x_{(i)}\right)^2 \right\} - \beta\mC^\alpha \ln(C)^2 R_{T},\\
\label{sn:alpha:beta}
\frac{\partial^2 l\left(\alpha, \beta \mid \bx \right)}{\partial \alpha \partial \beta} &=& -\sum_{i=1}^r\left\{ \left(1+R_i\right) x_{(i)}^\alpha \ln\left( x_{(i)}\right) \right\} - \mC^\alpha \ln(C) R_{T}\\
\label{second:deri:beta}
\frac{\partial^2 l\left(\alpha, \beta \mid \bx \right)}{\partial \beta^2} &=& \frac{-r}{\beta^2}
\end{eqnarray}

Now, using Equations \eqref{second:deri:alpha}-\eqref{second:deri:beta}, the Fisher's information matrix $\bI \left(\alpha, \beta\right)$ can be formed by 
\begin{eqnarray}\label{fisher:inf:matrix}
\bI \left(\alpha, \beta\right)=E\left[
          \begin{array}{cc}
           -\frac{\partial^2 l\left(\alpha, \beta \mid \bx \right)}{\partial \aa^2} & -\frac{\partial^2 l\left(\alpha, \beta \mid \bx \right)}{\partial \alpha \partial \beta}\\
           -\frac{\partial^2 l\left(\alpha, \beta \mid \bx \right)}{\partial \alpha \partial \beta} & -\frac{\partial^2 l\left(\alpha, \beta \mid \bx \right)}{\partial \beta^2}
          \end{array}
        \right].
\end{eqnarray}
It is well-known that (see \citet{Lawless}) the distribution of MLEs $\ll \ahat, \bhat\rr$ is a bivariate normal distribution with $$N\ll \ll \aa, \bb\rr , ~\bI^{-1}\ll \aa, \bb\rr \rr$$ where $\bI^{-1}\ll \aa, \bb\rr$ is the covariance matrix. Moreover, one can approximate the covariance matrix evaluated at $(\ahat,\bhat)$ by the following observed information matrix \begin{eqnarray}\label{obs:inf:matrix}
\bI \ll\ahat, \bhat \rr=\left[
          \begin{array}{cc}
           -\frac{\partial^2 l\left(\alpha, \beta \mid \bx \right)}{\partial \aa^2} & -\frac{\partial^2 l\left(\alpha, \beta \mid \bx \right)}{\partial \alpha \partial \beta}\\
           -\frac{\partial^2 l\left(\alpha, \beta \mid \bx \right)}{\partial \alpha \partial \beta} & -\frac{\partial^2 l\left(\alpha, \beta \mid \bx \right)}{\partial \beta^2}
          \end{array}
        \right]_{(\ahat,\bhat)}.
\end{eqnarray}

\subsection{Expectation-Maximization Algorithm}

The EM algorithm proposed by \citet{em1977} can be used to obtain the MLEs of the parameters $\alpha$ and $\beta$. It is known that the EM algorithm converges more reliably than NR.
Since Type-I PHCS can be considered as an incomplete data problem (see \citet{Ng2002}), it is possible to apply EM algorithm to obtain the MLEs of the parameters. Now, let us denote the incomplete (censored) data by $\bZ=\left(Z_1,Z_2,..., Z_r \right)$ where $Z_j=\left(Z_{j1}, Z_{j2}, ..., Z_{j R_j}  \right)$, $j=1,2,...,r$ such that $Z_j$ denotes the lifetimes of censored units at the time of $x_{(j)}$. Similarly, let $Z_{T}$ denotes the lifetimes of censored units at the time of $T$. Now, combining both the observed and censored data, one can obtain the complete data which is given by $\bW=\ll \bX, \bZ \rr$. The corresponding likelihood equation of the complete data can be obtained as follows:
\begin{eqnarray}\label{likelihood:complete}
L_W(\alpha, \beta|\bx) = \prod_{i=1}^{r}\left\{ f(\xi;\aa,\bb) \prod_{j=1}^{R_i} f(z_{i j};\aa,\bb) \right\} \prod_{j=1}^{R_T} f(z_{T j};\aa,\bb)
\end{eqnarray}
Therefore, the log-likelihood equation can be easily obtained by taking the natural logarithm of Equation \eqref{likelihood:complete} as follows:
\begin{eqnarray}\label{log-like:complete}
l_W(\alpha, \beta|\bx)&=&\ln\left(L_W(\alpha, \beta|\bx) \right)=\sum_{i=1}^{r}\ln\left(\alpha \beta \xi^{\alpha-1} {\rm exp}\left\{ -\beta \xi ^\alpha \right\} \right) + \sum_{i=1}^{r}\sum_{j=1}^{R_i}\ln \left(\alpha \beta z_{ij}^{\alpha -1}{\rm exp}\left\{ -\beta z_{ij}^{\alpha} \right\} \right) \no \\
&&+ \sum_{j=1}^{R_T}\ln \left(\alpha \beta z_{Tj}^{\alpha -1}{\rm exp}\left\{ -\beta z_{Tj}^{\alpha} \right\} \right) \no \\
&=& n \ln \alpha + n \ln \beta + (\alpha -1)\sum_{i=1}^{r}\ln\left(\xi\right)-\beta \sum_{i=1}^{r}\xi^\alpha + (\alpha -1)\sum_{i=1}^{r}\sum_{j=1}^{R_i}\ln \left(z_{ij}\right)-\beta\sum_{i=1}^{r}\sum_{j=1}^{R_i}z_{ij}^\alpha \no \\
&&+(\alpha -1)\sum_{j=1, r\neq m}^{R_T}\ln \left(z_{Tj}\right)-\beta\sum_{j=1, r\neq m}^{R_T}z_{Tj}^\alpha
\end{eqnarray}
Note that the last two terms of Equation\eqref {log-like:complete}, should be considered only for the Case II. Based on the complete sample, the MLEs of the parameters $\alpha$ and $\beta$ can be obtained by taking the derivatives of \eqref{log-like:complete} with respect to $\alpha$ and $\beta$ respectively and equating them to zero as follows:
\begin{eqnarray}\label{diff:alpha:comp}
\frac{\partial l_W(\alpha, \beta|\bx)}{\partial \alpha} &=& \frac{n}{\alpha}+ \sum_{i=1}^{r}\ln\left(\xi\right)-\beta \sum_{i=1}^{r} \xi ^\alpha \ln\left(\xi\right)+\sum_{i=1}^{r}\sum_{j=1}^{R_i}\ln \left(z_{ij}\right)-\beta \sum_{i=1}^{r}\sum_{j=1}^{R_i}z_{ij}^\alpha \ln \left(z_{ij}\right)\no \\
&&+\sum_{j=1, r\neq m}^{R_T}\ln \left(z_{Tj}\right)-\beta\sum_{j=1, r\neq m}^{R_T} z_{Tj}^\alpha \ln \left(z_{Tj}\right) = 0,\\
\label{diff:beta:comp}
\frac{\partial l_W(\alpha, \beta|\bx)}{\partial \beta} &=& \frac{n}{\beta}- \sum_{i=1}^{r}\xi^\alpha - \sum_{i=1}^{r}\sum_{j=1}^{R_i} z_{ij}^\alpha - \sum_{j=1, r\neq m}^{R_T} z_{Tj}^\alpha =0.
\end{eqnarray}
Now, the conditional expectation of the log-likehood equation of the complete data given the observations should be computed in the E-step of the algorithm. However, the following conditional expectations are necessary to be computed:
\begin{eqnarray}\label{cond:expect:alpha}
E\left(\frac{\partial l_W(\alpha, \beta|\bx)}{\partial \alpha} \Big|~ \xi, T \right) &=& \frac{n}{\alpha}+ \sum_{i=1}^{r}\ln\left(\xi\right)-\beta \sum_{i=1}^{r} \xi ^\alpha \ln\left(\xi\right)\no\\
&&+\sum_{i=1}^{r}\sum_{j=1}^{R_i} E\left[\ln \left(Z_{ij}\right) \left( 1- \beta Z_{ij}^\alpha \right) \Big| ~ Z_{ij} > \xi \right]\no\\
&&+\sum_{j=1, r\neq m}^{R_T}E\left[\ln \left(Z_{Tj}\right) \left(1-  \bb Z_{Tj}^\alpha \right) \Big| Z_{Tj} > T \right],
\end{eqnarray}
\begin{eqnarray}\label{cond:expect:beta}
E\left(\frac{\partial l_W(\beta, \beta ~|~ \bx)}{\partial \beta} \Big|  \xi, T \right) &=& \frac{n}{\beta}- \sum_{i=1}^{r}\xi^\alpha - \sum_{i=1}^{r}\sum_{j=1}^{R_i} E \left[Z_{ij}^\alpha \Big|  Z_{ij}>\xi\right]\no \\
&&- \sum_{j=1, r\neq m}^{R_T} E\left[Z_{Tj}^\alpha \Big| Z_{Tj}>T\right] .
\end{eqnarray}
In order to compute the expectations given above, making use of the theorem proved in \citet{Ng2002}, the conditional probability function of the censored data given the observed data can be obtained as follows:
\begin{eqnarray}\label{cond:pdf}
f(z_{i}|\mathcal{C}^{\ast},\alpha,\beta)&=&\frac{f(z_{i},\alpha,\beta)}{1-F(\mathcal{C}^{\ast},\alpha,\beta)}, Z_{i}>\mC^{\ast}
\end{eqnarray}
such that $\mC^{*} = x_{(i)}$ for $i = 1, 2, \ldots,
r$ and $\mathcal{C}^{*} = T$ for $i=T$.
Thus, the following expectations can be obtained
\begin{eqnarray}
\mE_{1}\left(\mC^*,\alpha, \beta \right)  &=& E\left[Z^\alpha \Big| Z > \mC^* \right] = \frac{1}{1-F\ll \mC^*,\alpha,\beta\rr} \int_{\mC^*}^{\infty}t^{\alpha}f(t) dt \no\\
&=& \frac{e^{-\bb\mC^{*\aa}}}{1-F\ll \mC^*,\alpha,\beta\rr}\ll 1+\bb\mC^{*\aa}\rr,
\end{eqnarray}
\begin{eqnarray}\label{Exp2}
\mE_{2}\left(\mC^*,\alpha, \beta \right) &=& E\left(\ln( Z)\left(1-\beta Z^\alpha \right) \Big| Z > \mC^* \right) \no\\
&=&\frac{1}{1-F\ll \mC^*,\alpha,\beta\rr} \int_{\mC^*}^{\infty}\ln( t)\left(1-\beta t^{\alpha} \right)  f(t)dt,
\end{eqnarray}
Since it is hard to obtain a closed form solution to Equation \eqref{Exp2}, the integral is approximated via Monte Carlo integration method in the simulation.
After updating the missing data with the expectations above in the E-step, the log-likelihood function is maximized in the M-step at the current state, say $\ahat_k$ and $\bhat_k$ being the estimators of $\alpha$ and $\beta$ and the following updating equations are computed:
\begin{eqnarray}
\ahat_{k+1} &=& n\left\{ -\sum_{i=1}^{r}\ln\left(\xi\right)+\bhat_{k+1} \sum_{i=1}^{r} \xi ^{\ahat_k} \ln\left(\xi\right)-\sum_{i=1}^r R_i\mE_2\left(\xi,\ahat_k, \bhat_{k+1}  \right)-R_T\mE_2\left(T,\ahat_k, \bhat_{k+1}  \right)
\right\}^{-1}\no\\
\label{EM:update}
\bhat_{k+1} &=& n\left\{ \sum_{i=1}^{r}\xi^{\ahat_k}+\sum_{i=1}^r R_i\mE_1\left(\xi,\ahat_k, \bhat_k  \right)+R_T\mE_1\left(T,\ahat_k, \bhat_k  \right)
\right\}^{-1}.
\end{eqnarray}
The EM estimates of $\left( \alpha, \beta \right)$ can be computed by an iterative procedure using Equation \eqref{EM:update} and the iterations can be terminated when $\abs{\ahat_{k+1} - \alpha_{k}} + \abs{ \bhat_{k+1} - \beta_{k}}<\epsilon$ where $\epsilon>0$ is a small real number. 

\subsection{Stochastic Expectation-Maximization Algorithm}
The computations in the E-step of EM algorithm is complex. Therefore, 
\citet{WeiTanner1990} proposed a Monte Carlo version of EM algorithm. However, the M-step of this algorithm may take so much time. \citet{EM1993} introduced a stochastic-EM (SEM) algorithm by considering a simulated values from the conditional distribution. \citet{Reza2018} used this algorithm successfully. In the SEM algorithm, firstly, one needs to generate $R_i$ number of samples of $z_{ij}$ where $i=1,2,...,r$ and $j=1,2,...,R_i$ using the following conditional CDF
\begin{eqnarray}\label{cond:cdf}
F\left( z_{ij}; \alpha, \beta | z_{ij} >\xi \right) &=& \frac{F\ll z_{ij};\alpha,\beta \rr-F\ll \xi;\alpha,\beta \rr}{1-F\ll \xi;\alpha,\beta \rr}, z_{ij}> \xi.
\end{eqnarray}
Now, using Equations \eqref{diff:alpha:comp} and \eqref{diff:beta:comp}, the estimators of $\ll \alpha, \beta\rr$ at the $k+1$ step of the algorithm can be obtained as follows:
\begin{eqnarray}\label{sem:beta}
\ahat_{k+1} &=& n\left[-\sum_{i=1}^{r}\ln\left(\xi\right)+\bhat_{k+1}\sum_{i=1}^{r} \xi ^{\ahat_k} \ln\left(\xi\right)-\sum_{i=1}^{r}\sum_{j=1}^{R_i}\ln \left(z_{ij}\right)\ll 1- \bhat_{k+1}z_{ij}^{\ahat_{k}} \rr 
\no \right.\\
 &&\left.-\sum_{j=1, r\neq m}^{R_T}\ln \left(z_{Tj}\right) \ll 1-  \bhat_{k+1} z_{Tj}^{\ahat_{k} }\rr \right]^{-1}\\
\label{sem:alpha}
\bhat_{k+1} &=& n\left[\sum_{i=1}^{r}\xi^{\ahat_{k}} + \sum_{i=1}^{r}\sum_{j=1}^{R_i} z_{ij}^{\ahat_{k}} + \sum_{j=1, r\neq m}^{R_T} z_{Tj}^{\ahat_{k}} \right]^{-1}.
\end{eqnarray}
Similarly, the iterations can be terminated when $\abs{\ahat_{k+1} - \alpha_{k}} + \abs{\bhat_{k+1} - \beta_{k}}<\epsilon$ where $\epsilon>0$ is a small real number.

\subsection{Fisher Information Matrix}
In this subsection, by making use of the idea of missing information principle proposed by \citet{Louis1982}, we can  obtain the observed Fisher
information matrix. \citet{Louis1982} suggested the following relation
\begin{eqnarray}\label{info:matrix}
\bI_{X}\ll\psi\rr = \bI_{W}\ll\psi\rr - \bI_{W \mid X}\ll\psi\rr
\end{eqnarray}
where $\psi=\ll\alpha, \beta\rr'$, $\bI_{X}\ll\psi\rr, \bI_{W}\ll\psi\rr$ and
$\bI_{W \mid X}\ll\psi\rr$ are the observed, complete and missing information matrices respectively. Now, the complete information matrix of a complete data set following the Weibull distribution can be obtained as
\begin{eqnarray}\label{comp:inf:mat}
\bI_{W}\ll\psi\rr &=& - E\ll\frac{\partial^2 \ln\mL}{\partial \psi^2}\rr = E\left[
          \begin{array}{cc}          \frac{n}{\alpha^{2}}+\beta \sum_{i=1}^n x_i^\alpha & \sum_{i=1}^n x_i^\alpha\ln x_i  \\
           \sum_{i=1}^n x_i^\alpha\ln x_i & \frac{n}{\beta^{2}}
          \end{array}
        \right]=\left[
          \begin{array}{cc}
           b_{11} & b_{12}\\
           b_{21} & b_{22}
          \end{array}
        \right]
\end{eqnarray}
where
\begin{eqnarray*}
b_{11} &=& \frac{n}{\alpha^2}+n\alpha\beta^2 \int_0^\infty \frac{x^{2\alpha-1}\ln(x)}{\rm{exp}(\beta x^\alpha)}dx\\
b_{12} &=& b_{21} = n\alpha\beta \int_0^\infty \frac{x^{2\alpha-1}\ln(x)}{\rm{exp}(\beta x^\alpha)}dx\\
b_{22} &=& \frac{n}{\beta^2}
\end{eqnarray*}
and $\ln \mL \ll \psi \rr = n\ln\alpha+n\ln\beta +(\alpha
-1)\sum_{i=1}^{n} x_i+\beta\sum_{i=1}^n x_i^\alpha$ is the corresponding log-likelihood equation.  Moreover, the missing
information matrix $\bI_{W \mid X}\ll\psi\rr$ is given by
\begin{equation}\label{missing:inf:mat}
\bI_{W \mid X}\ll\psi\rr = \sum_{i=1}^{r}R_{i}\bI^{(i)}_{W\mid X}\ll\psi\rr + R_{T}\bI^{*}_{W\mid X}\ll\psi\rr
\end{equation}
where $\bI^{(i)}_{W\mid X}\ll\psi\rr$ and $\bI^{*}_{W\mid X}\ll\psi\rr$ are the information matrices of a single observation from a truncated Weibull distribution from left at $\xi$ and $T$ respectively, such that
\begin{eqnarray*}
\bI^{(i)}_{W\mid X}\ll\psi\rr &=& -E \ll\frac{\partial^2 \ln\mL}{\partial \psi^2}\ln \left\{ f\ll z_{ij};\psi | z_{ij} > \xi \rr \right\} \rr.
\end{eqnarray*}
Now to calculate the missing information matrix $\bI^{(i)}_{W\mid X}\ll\psi\rr$, the conditional distribution given in Equation \eqref{cond:pdf} is used to obtain the following 
\begin{eqnarray*}
L_f&=&\ln\ll f(z_{ij}\mid z_{ij}>\xi) \rr = \ln(\aa) + \ln(\bb) + (\aa-1)\ln(z_{ij}) - \bb z_{ij}^\aa + \bb \xi^\aa.
\end{eqnarray*}
The second partial derivatives of $L_f$ are obtained as follows
\begin{eqnarray*}
\frac{\partial^2 L_f}{\partial \aa^2} &=& -\frac{1}{\aa^2} - \bb z_{ij}^\aa \ln(z_{ij})^2 + \bb \xi^\aa \ln(\xi)^2 \no \\
\frac{\partial^2 L_f}{\partial \aa \partial \bb} &=& -z_{ij}^\aa \ln(z_{ij}) + \xi^\aa \ln(\xi)\\
\frac{\partial^2 L_f}{\partial \bb^2} &=& -\frac{1}{\bb^2}.
\end{eqnarray*}
Now, in order to obtain the information matrices, the negative expected values of the quantities above are computed respectively as follows
\begin{eqnarray*}
E\ll -\frac{\partial^2 L_f}{\partial \aa^2}\rr &=& \frac{1}{\aa^2} + \bb \mE_{4}\ll\xi,\aa, \bb \rr - \bb \xi^\aa \ln(\xi)^2 \no \\
E\ll-\frac{\partial^2 L_f}{\partial \aa \partial \bb} \rr &=& \mE_{3}\ll\xi,\aa, \bb \rr - \xi^\aa \ln(\xi)\\
E\ll-\frac{\partial^2 L_f}{\partial \bb^2}\rr &=& \frac{1}{\bb^2}
\end{eqnarray*}
where 
\begin{eqnarray*}
\mE_{3}\ll\mC^*,\alpha, \beta \rr = E\ll Z^\alpha \ln(Z) \mid Z > \mC^* \rr &=& \frac{1}{1-F\ll \mC^*,\alpha,\beta\rr} \int_{\mC^*}^{\infty}t^{\alpha}\ln(t)f(t) dt \\
\mE_{4}\ll\mC^*,\alpha, \beta \rr = E\ll Z^\alpha \ln(Z)^2 \mid Z > \mC^* \rr &=& \frac{1}{1-F\ll \mC^*,\alpha,\beta\rr} \int_{\mC^*}^{\infty}t^{\alpha}\ln(t)^2f(t) dt .
\end{eqnarray*}
Using similar arguments, the information matrix $\bI^{*}_{W\mid X}\ll\psi\rr$ can also be computed easily. Then, using \eqref{info:matrix}--\eqref{comp:inf:mat}, the asymptotic variance-covariance matrix of $\widehat{\psi}$ can be computed by inverting the observed information matrix $\bI_{X}\ll\widehat{\psi}\rr$. Note that $\widehat{\psi}$ is computed using the NR estimates.

\section{Bayesian Estimation}\label{Bayes}

In this section, following \citet{Kundu2008}, we consider the Bayesian estimation for the parameters of the Weibull distribution under the assumption that the random variables $\aa$ and $\bb$ have independent gamma priors such that $\aa \sim Gamma(a,b)$ and
$\bb \sim Gamma(c, d)$. Therefore, the joint prior density of $\aa$ and $\bb$ can be written as
\begin{eqnarray*}\label{indep.prior}
\pi\ll \aa, \bb\rr \propto \aa^{a-1}\bb^{c-1}{\rm exp}\{ -(b\aa+d\bb)\},~a,b,c,d>0.
\end{eqnarray*}
Now, the posterior distribution of $\aa$ and $\bb$ can be obtained as follows
\begin{eqnarray}\label{indep:posterior}
\pi\ll \alpha, \beta \mid \bx\rr &=& \frac{L\ll \alpha, \beta \mid \bx\rr \pi\ll \alpha, \beta\rr}{\int_0^\infty \int_0^\infty L\ll \alpha, \beta \mid \bx\rr \pi\ll \alpha, \beta\rr d\aa d\bb}\no \\
&=& \frac{\ll\prod_{i=1}^r \xi^{\aa-1}\rr \bb^{c+r-1} \aa^{a+r-1}}{\Gamma(c+r)\Psi(a,c,\bx)} {\rm exp}\left\{d-b\aa+ \sum_{i=1}^r(1+R_i)\xi^{\aa}+C^\aa R_T \right\}
\end{eqnarray}
where 
\begin{eqnarray*}
\Psi(a,c,\bx)=\int_0^\infty\frac{\aa^{a+r-1}{\rm exp}\left\{ -b\aa\right\} \ll\prod_{i=1}^r \xi^{\aa-1}\rr}{\left[d+ \sum_{i=1}^r(1+R_i)\xi^{\aa}+C^\aa R_T \right]^{a+c+r}}d\aa.
\end{eqnarray*}

In this paper, three different loss functions are considered. One of them is the most commonly used squared error loss function (SEL) which is defined as follows:
\begin{eqnarray*}
L_S\ll\widehat{t}(\psi), t(\psi)\rr = \ll\widehat{t}(\psi)-t(\psi)\rr^2
\end{eqnarray*}
where $\widehat{t}(\psi)$ is an estimator of $t(\psi)$. SEL is a symmetric loss function which gives equal weights to both underestimation and overestimation.
Secondly, linear-exponential (LINEX) loss function which is a useful asymmetric loss function introduced by
\citet{linex} as follows
\begin{eqnarray*}\label{linex}
L_L\ll\widehat{t}(\psi), t(\psi)\rr =e^{\nu\ll\widehat{t}(\psi)-t(\psi)\rr}-\nu\ll\widehat{t}(\psi)-t(\psi)\rr-1, \;\nu\neq0.
\end{eqnarray*}
The LINEX loss function is a convex function whose shape is determined by the value of $\nu$. The negative (positive) value of $\nu$ gives more weight to overestimation (underestimation) and its magnitude reflects the degree of asymmetry. It is seen that, for $\nu=1$, the function is quite asymmetric with overestimation being costlier than underestimation. If $\nu<0$, it rises almost exponentially when the estimation error $\widehat{t}(\psi)-t(\psi)<0$ and almost linearly if $\widehat{t}(\psi)-t(\psi)>0$.  For small values of $|\nu|$, the LINEX loss function is almost symmetric and not far from squared error loss function.

Under the SEL function, the Bayes estimators of $\aa$ and $\bb$ which are the expected values of the corresponding posterior distributions are computed respectively as follows
\begin{eqnarray}
\ahat_S= E\ll\pi\ll \alpha \mid \bx\rr \rr=\frac{\Psi(a+1,c-1,\bx)}{\Psi(a,c,\bx)}
\end{eqnarray}
and
\begin{eqnarray}
\bhat_S=E\ll\pi\ll \bb \mid \bx\rr \rr= (a+c+r)\frac{\Psi(a,c+1,\bx)}{\Psi(a,c,\bx)}.
\end{eqnarray}
Since the Bayes estimators given above includes the complicated integral function $\Psi(a,c+1,\bx)$ we also consider using the Bayes estimate of $t(\psi)$ under the LINEX loss function is given by
\begin{eqnarray*}\label{llbays}
\widehat{t}_{L}(\psi)=-\frac{1}{\nu}\ln\left[E_{t}\ll e^{-\nu t(\psi)}\mid\bm{x} \rr \right]~=~ -\frac{1}{\nu}\ln\left[\int_{0}^{\infty}\int_{0}^{\infty}e^{-\nu t(\psi)}\pi(\aa, \beta \mid \bm{x})d\aa \, d\beta\right].
\end{eqnarray*}
Finally, the general entropy loss (GEL) function is also considered and it is given by
\begin{eqnarray*}\label{gen}
L_{GEL}\ll\widehat{t}(\psi), t(\psi)\rr
=\ll\frac{\widehat{t}(\psi)}{t(\psi)}\rr^{\kappa}-\kappa\ln\ll\frac{\widehat{t}(\psi)}{t(\psi)}\rr -1,\: \kappa\neq 0.
\end{eqnarray*}
where $\kappa$ is the shape parameter showing the departure from symmetry. When $\kappa>0$, the overestimation is considered to be more serious than underestimation and for $\kappa<0$ vice versa. The Bayes estimator under GEL function is given by
\begin{eqnarray*}
\widehat{t}_{GEL}(\psi)=\left[E_{t}\ll  t(\psi)^{-\kappa}\mid\bm{x} \rr \right]^{-1/\kappa}~=~ \left[\int_{0}^{\infty}\int_{0}^{\infty}t(\psi)^{-\kappa}\pi(\aa, \beta \mid \bm{x})d\aa \, d\beta\right]^{-1/\kappa}.
\end{eqnarray*}

\subsection{Tierney-Kadane Approximation}
In this subsection, the approximation method of \citet{Tierney} is used to obtain the approximate Bayes estimators under SEL, LINEX and GEL loss functions. Now, we consider the following functions
\begin{eqnarray}\label{TK:Delta:func}
\Delta(\aa,\beta) &=& \frac{1}{n}\ln[L(\aa,\bb\mid\bm{x})\pi(\aa,\beta)],\\ \label{TK:Delta*:func}
\Delta^{*}(\aa,\bb) &=& \frac{1}{n}\ln[L(\aa,\bb\mid\bm{x})\pi(\aa,\bb)t(\psi)].
\end{eqnarray}
Now assume that $(\atil_{\Delta}, \btil_{\Delta})$ and  $(\atil_{\Delta^{*}}, \btil_{\Delta^{*}})$ respectively maximize the functions $\Delta(\aa,\bb)$ and $\Delta^{*}(\aa,\bb)$. Then the approximation method of \citet{Tierney} is given by
\begin{eqnarray*}\label{TK-est:SEL}
\widetilde{t}_{\rm SEL}(\aa,\bb)=\sqrt{\dfrac{\vert \Sigma^{\ast}\vert}{\vert \Sigma\vert}}exp\left[ n \ll \Delta_1^{*}\ll\atil_{\Delta^{*}},\btil_{\Delta^{*}}\rr - \Delta\ll\atil_{\Delta},\btil_{\Delta}\rr  \rr\right]
\end{eqnarray*}
where $|\Sigma|$ and $|\Sigma^{*}|$ are the negative of inverses the second derivative matrices of $\Delta\ll\aa,\bb\rr$ and $\Delta_1^* \ll\aa,\bb\rr$ respectively obtained at $(\atil_{\Delta}, \btil_{\Delta})$ and  $(\atil_{\Delta^{*}}, \btil_{\Delta^{*}})$. The function $\Delta(\aa,\bb)$  can be easily obtained by using the Equation \eqref{TK:Delta:func} as follows
\begin{eqnarray}\label{delta}
\Delta(\aa,\beta) &=& \frac{1}{n}\left[\ln(M)+(\aa-1)\sum_{i=1}^r\ln(\xi) -\bb \ll d+b\aa+ \sum_{i=1}^r(1+R_i)\xi^{\aa}+C^\aa R_T \rr \no \right.\\
 &&\left.+ (a+c+r-1)\ln(\bb)+(a+r-1)\ln(\aa)\right]
\end{eqnarray}
where $M=\frac{d^c b^a}{\Gamma(c)\Gamma(a)}$. Now, differentiating Equation \eqref{delta} with respect to $\aa$ and $\bb$ solving for these parameters, one gets the following equations
\begin{eqnarray*}
\atil_{\Delta} &=& (a+r-1)\left[\bb \ll b+ \sum_{i=1}^r(1+R_i)\xi^{\aa}+C^\aa R_T \rr-\sum_{i=1}^r\ln(\xi) \right]^{-1},\\
\btil_{\Delta} &=& (a+c+r-1)\left[ \sum_{i=1}^r(1+R_i)\xi^{\aa}+C^\aa R_T +d+b\aa\right]^{-1}.
\end{eqnarray*}
Since it is easy to obtain the second derivatives and the related Hessian matrices, we skip this part. Thus under the SEL function, the approximate Bayes estimators are computed by
\begin{eqnarray*}\label{TK:alpha:SEL}
\atil_{\rm SEL} &=& \sqrt{\dfrac{\vert \Sigma^{\ast}\vert}{\vert \Sigma\vert}}exp\left[ n \ll \Delta_{1\aa}^{*}\ll\atil_{\Delta^{*}},\btil_{\Delta^{*}}\rr - \Delta\ll\atil_{\Delta},\btil_{\Delta}\rr  \rr\right], \\ \label{TK:beta:SEL}
\btil_{\rm SEL} &=& \sqrt{\dfrac{\vert \Sigma^{\ast}\vert}{\vert \Sigma\vert}}exp\left[ n \ll \Delta_{1\bb}^{*}\ll\atil_{\Delta^{*}},\btil_{\Delta^{*}}\rr - \Delta\ll\atil_{\Delta},\btil_{\Delta}\rr  \rr\right]
\end{eqnarray*}
where $\Delta_{1\aa}^{*}\ll \aa,\bb\rr = \Delta\ll \aa,\bb\rr+\frac{1}{n}\ln(\aa)$ for $t(\aa,\bb) = \aa$ and $\Delta_{1\bb}^{*}\ll \aa,\bb\rr = \Delta\ll \aa,\bb\rr+\frac{1}{n}\ln(\bb)$ for $t(\aa,\bb) = \bb$.

One can also compute the Bayes estimators under the LINEX loss and get
\begin{eqnarray*}\label{TK-est:LINEX}
\widetilde{t}_{\rm LINEX}(\aa,\bb)=\sqrt{\dfrac{\vert \Sigma^{\ast}\vert}{\vert \Sigma\vert}}exp\left[n\left\{ \Delta_2^{*}\ll \atil_{\Delta^{*}},\btil_{\Delta^{*}}\rr -\Delta\ll \atil_{\Delta},\btil_{\Delta}\rr\right\}\right].
\end{eqnarray*}
Letting $t(\aa,\bb)=e^{-\nu \aa}$, one gets  $\Delta_{2\aa}^{*}\ll \aa,\bb\rr=\Delta\ll \aa,\bb\rr-\frac{1}{n}\nu\aa$ and letting $t(\aa,\bb)=e^{-\nu \bb}$, $\Delta_{2\bb}^{*}\ll \aa,\bb\rr= \Delta\ll \aa,\bb\rr-\frac{1}{n}\nu\bb$. Thus, approximate Bayes estimators under LINEX function are computed as

\begin{eqnarray*}\label{TK:alpha:LIN}
\atil_{\rm LINEX} &=& -\frac{1}{\nu}\ln\ll\sqrt{\dfrac{\vert \Sigma^{\ast}\vert}{\vert \Sigma\vert}}exp\left[ n \ll \Delta_{2\aa}^{*}\ll\atil_{\Delta^{*}},\btil_{\Delta^{*}}\rr - \Delta\ll\atil_{\Delta},\btil_{\Delta}\rr  \rr\right]\rr, \\ \label{TK:beta:LIN}
\btil_{\rm LINEX} &=& -\frac{1}{\nu}\ln\ll\sqrt{\dfrac{\vert \Sigma^{\ast}\vert}{\vert \Sigma\vert}}exp\left[ n \ll \Delta_{2\bb}^{*}\ll\atil_{\Delta^{*}},\btil_{\Delta^{*}}\rr - \Delta\ll\atil_{\Delta},\btil_{\Delta}\rr  \rr\right]\rr.
\end{eqnarray*}

Finally, letting $t(\aa,\bb)=\aa^{-\kappa}$, one gets  $\Delta_{3\aa}^{*}\ll \aa,\bb\rr=\Delta\ll \aa,\bb\rr-\frac{\kappa}{n}\ln(\aa)$ and letting $t(\aa,\bb)=\bb^{-\kappa}$, $\Delta_{3\bb}^{*}\ll \aa,\bb\rr= \Delta\ll \aa,\bb\rr-\frac{\kappa}{n}\ln(\bb)$. Thus, approximate Bayes estimators under GEL function are obtained by
\begin{eqnarray*}\label{TK:alpha:GEL}
\atil_{\rm GEL} &=& \ll\sqrt{\dfrac{\vert \Sigma^{\ast}\vert}{\vert \Sigma\vert}}exp\left[ n \ll \Delta_{3\aa}^{*}\ll\atil_{\Delta^{*}},\btil_{\Delta^{*}}\rr - \Delta\ll\atil_{\Delta},\btil_{\Delta}\rr  \rr\right]\rr^{-1/\kappa}, \\ \label{TK:beta:GEL}
\btil_{\rm GEL} &=& \ll\sqrt{\dfrac{\vert \Sigma^{\ast}\vert}{\vert \Sigma\vert}}exp\left[ n \ll \Delta_{3\bb}^{*}\ll\atil_{\Delta^{*}},\btil_{\Delta^{*}}\rr - \Delta\ll\atil_{\Delta},\btil_{\Delta}\rr  \rr\right]\rr^{-1/\kappa}.
\end{eqnarray*}

\subsection{MCMC Method}
Metropolis--Hastings (MH) algorithm, a method for generating random samples from the posterior distribution using a proposal density, is considered in this subsection. A symmetric proposal density of type $q(\theta' | \theta ) = q(\theta | \theta')$ may be considered generally, where $\theta$ is the parameter vector of the distribution considered. Following \citet{Dey2016b}, we consider a bivariate normal distribution as the proposal density such that $q(\theta' | \theta ) = N(\theta | \bV_\theta)$ where $\bV_\theta$ is the covariance matrix and $\theta=(\aa,\bb)$. Although, the bivariate normal distribution may generate negative observations, the domain of both shape and scale parameters of Weibull distribution is positive. Therefore, the following steps of MH algorithm is used to generate MCMC sample from the posterior density given by \eqref{indep:posterior}
\begin{itemize}
    \item[(1)] Set the initial parameter values as $\theta=\theta_0.$
    \item[(2)] For $j=1,2,...,N$, repeat the following steps:
    \begin{itemize}
        \item[(i)] Set $\theta=\theta_{j-1}$
        \item[(ii)] Generate new parameters $\lambda$ from bivariate normal $N_2\ll \ln(\theta),\bV_\theta \rr$ \item[(iii)] Compute $\theta_{new}={\rm exp}(\lambda)$
        \item[(iv)] Calculate $\gamma={\rm min}\ll 1, \frac{\pi\ll \theta_{new} \mid \bx\rr \theta_{new}}{\pi\ll \theta \mid \bx\rr \theta} \rr$
        \item[(v)] Set $\theta_j=\theta_{new}$ with probability $\lambda$, otherwise $\theta_j=\theta.$
    \end{itemize}
\end{itemize}
After generating the MCMC sample, some of the initial samples, say $N_0$, can be discarded as burn-in process and the estimations can be computed via the remaining ones $(M=N-N_0)$ under SEL, LINEX and GEL loss functions as follows
\begin{eqnarray*}
\widehat{t}_{\rm SEL}(\psi)&=&\frac{1}{M}\sum_{i=1}^Mt(\psi_i), \\
\widehat{t}_{\rm LINEX}(\psi) &=& -\frac{1}{\nu}\ln\ll \frac{1}{M}\sum_{i=1}^M exp\ll -\nu t(\psi_i)\rr\rr, \\
\widehat{t}_{\rm GEL}(\psi)&=&\ll \frac{1}{M}\sum_{i=1}^M \ll t(\psi_i)^{-\kappa}\rr \rr^{-1/\kappa}.
\end{eqnarray*}
The main advantage of MCMC method over Tierney--Kadane method is that the MCMC samples can also be used to compute highest posterior density (HPD) intervals. \citet{chen-shao} proposed a method to compute the HPD intervals using MCMC samples. This method has been used in the literature extensively. Now, consider the posterior density $\pi(\theta | \bx)$. Assume that the pth quantile of the distribution is given by $\theta^{(p)}={\rm inf}\left\{\theta:\Pi(\theta|\bx) \geq p; 0<p<1 \right\}$ where $\Pi(\theta|\bx)$ denotes the posterior distribution function of $\theta$. Now, for a given $\theta^*$, a simulation consistent estimator of $\Pi(\theta^* | \bx)$ can be computed as
\begin{eqnarray*}
\Pi(\theta^*|\bx)= \frac{1}{M}\sum_{i=1}^M I(\theta \leq \theta^*)
\end{eqnarray*}
where $I(\theta \leq \theta^*)$ is an indicator function. Then, the estimate of $\Pi(\theta^*|\bx)$ is given as
\begin{eqnarray*}
\widehat{\Pi}(\theta^*|\bx) = 
\begin{cases}
      0 & \text{if} ~\theta^*<\theta_{(N_0)}\\
      \sum_{j=N_0}^i \gamma_j & \text{if} ~\theta_{(i)}<\theta^*<\theta_{(i+1)}\\
      1 & \text{if} ~\theta_{(M)}
    \end{cases}  
\end{eqnarray*}
where $\gamma_j=1/M$ and $\theta_{(j)}$ is the jth ordered value of $\theta_{j}$.
$\theta^{(p)}$ can be approximated by the following
\begin{eqnarray*}
\theta^{(p)}=\begin{cases}
      \theta_{(N_0)} & \text{if}~ p=0\\
      \theta_{(j)} & \text{if}~\sum_{j=N_0}^{i-1} \gamma_j<p<\sum_{j=N_0}^i \gamma_j
\end{cases}
\end{eqnarray*}
Now, one can construct the $100(1-p)\%$ confidence intervals where $0<p<1$ as $\ll\widehat{\theta}^{j/s}, \widehat{\theta}^{(j+[(1-p)s])/s} \rr$, $j=1, 2,..., s-[(1-p)s]$ such that $[v]$ denotes the greatest integer less than or equal to $v$. At the end, the HPD credible interval of $\theta$ is the one having the shortest length.

\section{Shrinkage Estimation}\label{Shrinkage}

Prior information on the parameters in a statistical model generally leads to an improved inference procedure in problems of statistical inference. Restricted models arise from the incorporation of the known prior information in the model in the form of a constraint. The estimators obtained from restricted (unrestricted) model is known as the restricted (unrestricted) estimators. 
The results of an analysis of the restricted and unrestricted models can be weighted against loss of efficiency and validity of the constraints in deciding a choice between these two extreme inference methods, when a full confidence may not be in the prior information, see \citet{ahmed1990}.

\citet{bancroft} was the first to consider a pre-test procedure when there is doubt that the prior information is not certain (uncertain prior information). After the pioneering study of \citet{bancroft}, pre-test estimators has gained much attention. \citet{tompson} defined an efficient shrinkage estimator. Following \citet{tompson}, shrinkage estimation of the Weibull parameters has been discussed by a number of authors, including \citet{SB1978}, \citet{Pandey1983}, \citet{PS1993} and \citet{SS2000}. We also refer to the following book and papers among others: \citet{JB}, \citet{SalehKibria1993}, \citet{Saleh2006}, \citet{kibria2010}.

Now suppose that there is an uncertain prior information in the form of $\theta=\theta_0$ where $\theta$ is the parameter of a distribution of interest. Our aim is to estimate $\theta$ using a pre-test estimation strategy and this prior information. Therefore, we consider the following hypothesis to check the validity of this information
\begin{eqnarray*}
H_0: \theta=\theta_0\\ 
H_0: \theta \neq \theta_0\no
\end{eqnarray*}
It is known that under $H_0$, the asymptotic distribution of $\sqrt{D}(\widehat{\theta}-\theta_0
)$ is normal with $N(0, \sigma_{\widehat{\theta}}^2)$ and the related test statistics can be defined as follows
\begin{eqnarray*}
W_D = \left(\frac{\sqrt{D}(\widehat{\theta}-\theta_0
)}{\sigma_{\widehat{\theta}}^2} \right)^2.
\end{eqnarray*}
One can reject the null hypothesis when $W_D>\chi^2_{1}(\lambda)$ based on the distribution of $W_D$ where $\lambda$ can be treated as the degree of trust in the prior information about the parameter such that $\theta=\theta_0$. Thus, the shrinkage pre-test estimator (SPT) can be defined as
\begin{eqnarray*}
\widehat{\theta}_{\rm SPT} = \lambda\theta_0 +(1-\lambda)\widehat{\theta} I\left(W_D<\chi^2_{1}(\lambda)\right)
\end{eqnarray*}
where $I(A)$ is the indicator of the set $A$.


\section{Monte Carlo Simulation Experiments}\label{sim}

In this section, we conduct a simulation study to illustrate the performance of the different estimation techniques discussed in this paper by considering $(n, m) = (30, 15), (30, 20), (50, 25), (50, 40)$, different values of predetermined time $T = 0.21, 1.15$ such that these values correspond to  the sample quantiles for the probabilities $0.5$ and $0.8$ respectively, and the real values of the parameters are chosen as $\aa = 0.5$ and $\bb = 1.5$ in all cases. The following three schemes are considered in the simulation
\begin{itemize}
    \item 
    Scheme 1: $R = (0^{m-1},n-m)$
    \item
    Scheme 2: $R = (n-m, 0^{m-1})$
    \item
    Scheme 3: $R = (2^5, 0^{m-6}, n-m-10)$
\end{itemize}
It is noted that Scheme 1 is the Type-II censoring such that $n-m$ units are removed from the experiment at the time of the m−th failure, in Scheme 2, $n-m$ units are removed at the time of the first failure. However, in Scheme 3,  a progressive Type-II censoring scheme allowing different numbers of censoring within the experiment is considered.
The progressive type II censored data from Weibull distribution is generated using
algorithm proposed by \citet{Bala2000}.
The maximum likelihood estimators of $\aa$ and $\bb$ are obtained using NR, EM and SEM algorithms.
In computing the Bayes estimates, two different priors are used such as the non-informative priors as $a=b=c=d=0$ and the informative priors where we assume that we have past samples from Weibull$(\aa,\bb)$ distribution, say $K$ samples and their corresponding MLEs as $\ll\widehat{\aa}_j,\widehat{\bb}_j\rr$, $j=1, 2, ..., K$. Now, equating the sample means and variances of these values to the means and variances of gamma priors respectively and solving the equations for $K=1000$, and $n=30$ being the sample size of past samples, we obtain the following informative prior values, $a=43.77, b=83.45, c=24.24, d=15.47$. 

Bayes estimates are computed under SEL, LINEX, GEL loss functions. Notice that for the LINEX loss function, we considered two values of $\nu$ as $\nu = -0.5, 0.5$ giving more weight to underestimation and overestimation respectively. Similarly, two choices of $\kappa$ such as $\kappa = -0.5, 0.5$ are taken into account under GEL function. Moreover, $6000$ MCMC samples are generated and MCMC estimations are computed under the listed loss function and respective parameter values. The first $1000$ MCMC samples are considered as a burn-in sample so that the average values and MSEs are computed via the remaining 5000 samples for each replicate in the simulation. 

For the shrinkage estimators, the test statistic $W_D$ is calculated and then shrinkage pre-test (SPT) estimators are obtained. The distribution of the test statistic $W_D$ is computed under the null hypothesis, that is, $H_0:\theta = \theta_0$. Moreover, we take $\lambda = 0.5$ giving equal weight to both restricted and unrestricted estimators and the type one test error is set to $0.05$ in testing the hypothesis, prior values of the parameters are taken as $\aa_0=0.7, \bb_0=1.7$ for practical purposes. The MLE shrinkage pre-test estimators are obtained using NR algorithm and also the Bayes estimator with T-K method under different loss functions.

Totally, $5000$ repetitions are carried out and average values (Avg), mean squared errors (MSE), confidence/ credible interval lengths (IL) and coverage probabilities (CP) are obtained for the purpose of comparison. MSEs of the estimators are computed as follows
\begin{eqnarray*}
{\rm MSE}\ll\widehat{\theta}\rr=\frac{1}{5000}\sum_{i=1}^{5000}\ll\widehat{\theta}_i-\theta \rr^2
\end{eqnarray*}
where $\widehat{\theta}_i$ is NR, EM, SEM, SPT estimators and Bayes estimators under SEL loss function in the ith replication. However, the MSEs of Bayes estimators under LINEX and GEL loss functions are computed respectively by
\begin{eqnarray*}
{\rm MSE_{LINEX}}\ll\widehat{\theta}\rr&=&\frac{1}{5000}\sum_{i=1}^{5000}\ll e^{\nu\ll\widehat{\theta}_i-\theta\rr}-\nu\ll\widehat{\theta}_i-\theta\rr-1\rr, \\
{\rm MSE_{GEL}}\ll \widehat{\theta}\rr&=&\frac{1}{5000}\sum_{i=1}^{5000}\ll\ll\frac{\widehat{\theta}_i}{\theta}\rr^{\kappa}-\kappa\ln\ll\frac{\widehat{\theta}_i}{\theta}\rr -1\rr.
\end{eqnarray*}
All of the computations are performed using the R Statistical Program \citep{R2018}. All the results are presented in Tables \ref{NR:n30}--\ref{CI:CP:n50}.

Based on Tables \ref{NR:n30} and \ref{NR:n50}, we can conclude that EM estimates are quiet preferable to the SEM and NR method for all schemes and $T$s. Both MSEs and Avgs for EM estimates are the smallest. However, NR method beats SEM in terms of lower Avgs and MSEs. We also observe that as $m$ increase, the values of MSEs and Avgs decrease, generally.

We reported the results of Bayes estimates based on TK and MCMC methods in Tables \ref{TK:n30m15}--\ref{MCMC:n50m40}. From these tables, it is evident that all the Bayes estimates based on informative priors have very small MSEs compared to the MLEs. We also see that the Bayes estimates based on informative priors are better than those that are based on non-informative priors in all schemes and $(T,n,m)$s. However, EM estimates are better than non-informative Bayes estimates based on SEL in terms of MSE when $n=30$ and in terms of both MSE and Avg when $n=50$.  So we can conclude that Bayes estimates even with non informative priors are preferable to the MLE except for EM estimates, for all schemes and $T$s. When we compare MSEs of T-K and MCMC methods, we observed that they are generally close to each other. However, T-K is better in some of the cases and vice versa in some others.  However, the MCMC has the advantage of construction the credible intervals. Thus, we can say that MCMC is preferable since it gives more information.

The performances of SPT estimators are given in Tables \ref{SPT:n30}--\ref{SPT:n50}. According to these tables, we can say that SPT estimators based on informative T-K method have better performance than SPT based on NR methods in the sense of both MSE and Avg, generally. Moreover, SPT with T-K method based on GEL function seems to have the least MSE values among others. SPT estimator based on NR method has smaller MSE values than NR estimator when we consider the parameter $\bb$, and both methods have closer MSE values for the parameter $\aa$.

Finally, we have summarized the confidence intervals and coverage probabilities in Tables \ref{CI:CP:n30}--\ref{CI:CP:n50} for $n=30$ and $n=50$, respectively. It is observed that when we use non-informative priors the estimated CPs are smaller than the nominal CPs. Moreover, the expected ILs of non-informative methods are less than that of NR method. However, the estimated CPs of NR are slightly more than the non-informative method. Further, we observe that the CIs based on informative priors are better than the ones based on the non-informative priors and the once based on NR, in terms of having smaller ILs but higher CPs. We can also find that when we increase $T$, it is seen that the ILs are generally decreasing. Finally, it is observed that for $n=50$, the estimated CIs are getting better in terms of having smaller ILs and higher CPs as compared to $n=30$.

\section{Real Data Example}\label{RealData}

We consider a data set reported by \citet{real_data} representing the strength measured in GigaPAscal (GPA) for single carbon fibres, and impregnated 1000-carbon fibre tows. Single fibres were tested under tension at gauge lengths of 10 mm. This data was analysed by \citet{Asgar} considering a hybrid censoring scheme for the Weibull distribution. Following \citet{Asgar}, we analyzed this data set using two-parameter Weibull distribution after subtracting 1.75. The authors recorded that the validity of the Weibull model based on the Kolmogorov--Smirnov (K--S) test is full-filled, namely, K–S = 0.072 and p-value = 0.885.

To compute the Bayes estimates, since we have no prior information about the unknown parameters, we assume the non-informative priors by setting $a=b=c=d=0$. Taking $m=40$ and $T=2$, we used the following schemes
\begin{itemize}
    \item 
    Scheme 1: $R = (0^{39},23)$
    \item
    Scheme 2: $R = (23, 0^{39})$
    \item
    Scheme 3: $R = (2, 0^{10},2^3,0^{10}, 2^3,0^{10},3^3)$
\end{itemize}
In SPT estimates, since we don't have any prior information about parameters, we use the Bayes estimates as a an estimated prior information. Then we substitute them in the SPT formulae as $\widehat{\theta}_{\rm SPT} = \lambda\theta_0 +(1-\lambda)\widehat{\theta}_{Bayes} I\left(W_D<\chi^2_{1}(\lambda)\right)$ by setting $\lambda=0.5$  and $\alpha=0.05$.

All the estimation methods considered in this paper are applied to this data and the estimated parameter values are reported in Table \ref{app:est}. We observe that the estimated values of $\aa$ and $\bb$ based on Bayes method are closer to each other while these values are a bit different than each other based on MLE methods. Further, it can be seen that the Bayes estimates based on non-informative priors for all loss functions are close to the NR estimates as compared to EM and SEM estimates.
Moreover, asymptotic confidence intervals of NR method and HPD intervals of MCMC method are given in Table \ref{app:CI}. According to this table, we can say that MCMC confidence intervals are mostly wider than the ones obtained via NR. However, in simulation study, we observed that the MCMC confidence intervals are preferable to NR confidence intervals.

\section{Conclusive Remarks}\label{conclusion}
In this paper, we discussed the estimation of parameters of Weibull distribution under Type-I progressively hybrid censoring scheme using both classical and Bayesian strategies. Namely, MLE is obtained using NR, EM and SEM algorithms and  Bayesian estimators are computed via T-K approximation and MCMC method under SEL, LINEX and GEL loss functions. We have also proposed the shrinkage preliminary test estimators based on NR and T-K with informative priors using equal weights on the prior information and the sample information. A real data application and extensive Monte Carlo simulations have been considered to compare the estimators in terms of MSE and Avg and also we compared the lengths of CIs and CPs. According to the results, EM algorithm beats the other ML estimates. However, we observed that both the T-K and MCMC methods perform quite closely. Finally, we found out that shrinkage preliminary test estimates have satisfactory  performances in the presence of having proper prior information. 
\vskip 1cm
\noindent\textbf{Acknowledgements.} This paper was written while Yasin Asar visited McMaster University and he was supported by The Scientific and Technological Research Council of Turkey (TUBITAK), BIDEB-2219 Postdoctoral Research Program, Project No: 1059B191700537.

\newpage
\begin{table}[ht] 
\centering
\caption{Average values and the corresponding MSEs of the estimators NR, EM and SEM when $n=30$.}
\label{NR:n30}
\begin{tabular}{ccccccccc}
\multicolumn{9}{c}{$m=15$}\\
  \hline
 &&& \multicolumn{2}{c}{NR} & \multicolumn{2}{c}{EM} & \multicolumn{2}{c}{SEM} \\ 
 \cmidrule{4-9} 
 $T$& $R$&&$\aa$&$\bb$&$\aa$&$\bb$&$\aa$&$\bb$ \\ 
 \hline
     & 1 & Avg & 0.5441 & 1.7863 & 0.5220 & 1.5999 & 0.5602 & 1.8180 \\ 
       & & MSE & 0.0218 & 0.6901 & 0.0099 & 0.1970 & 0.0251 & 0.6874 \\ 
0.21 & 2 &  Avg & 0.5322 & 1.7720 & 0.5216 & 1.5953 & 0.5679 & 1.8793 \\ 
       & & MSE & 0.0253 & 0.7994 & 0.0125 & 0.1912 & 0.0332 & 0.8777 \\ 
     & 3 & Avg & 0.5329 & 1.7287 & 0.5225 & 1.5935 & 0.5569 & 1.8119 \\ 
       & & MSE & 0.0210 & 0.5954 & 0.0106 & 0.2043 & 0.0259 & 0.6829 \\
  \hline
     & 1 & Avg & 0.5544 & 1.8471 & 0.5253 & 1.6365 & 0.5620 & 1.8506 \\ 
       & & MSE & 0.0216 & 0.6985 & 0.0107 & 0.2093 & 0.0243 & 0.7161 \\ 
1.15 & 2 & Avg & 0.5244 & 1.6475 & 0.5207 & 1.5987 & 0.5430 & 1.7002 \\ 
       & & MSE & 0.0151 & 0.3337 & 0.0133 & 0.1854 & 0.0208 & 0.4190 \\ 
     & 3 & Avg & 0.5472 & 1.7615 & 0.5343 & 1.6566 & 0.5563 & 1.7659 \\ 
       & & MSE & 0.0176 & 0.5286 & 0.0128 & 0.2769 & 0.0216 & 0.5452 \\ \hline 
\end{tabular}
\begin{tabular}{ccccccccc}
\multicolumn{9}{c}{$m=20$}\\
  \hline
     & 1 & Avg & 0.5380 & 1.6806 & 0.5254 & 1.5823 & 0.5586 & 1.7527 \\
       & & MSE & 0.0209 & 0.4365 & 0.0109 & 0.1730 & 0.0253 & 0.5087 \\ 
0.21 & 2 & Avg & 0.5330 & 1.7445 & 0.5241 & 1.6172 & 0.5591 & 1.8429 \\ 
       & & MSE & 0.0232 & 0.6426 & 0.0127 & 0.2177 & 0.0289 & 0.7446 \\ 
     & 3 & Avg & 0.5333 & 1.7409 & 0.5217 & 1.5959 & 0.5573 & 1.8283 \\ 
       & & MSE & 0.0221 & 0.6361 & 0.0106 & 0.2033 & 0.0267 & 0.7517 \\ 
       \hline 
     & 1 & Avg & 0.5442 & 1.7104 & 0.5306 & 1.6317 & 0.5537 & 1.7362 \\ 
       & & MSE & 0.0154 & 0.3554 & 0.0107 & 0.1921 & 0.0192 & 0.4150 \\
1.15 & 2 & Avg & 0.5283 & 1.6088 & 0.5262 & 1.5838 & 0.5433 & 1.6453 \\ 
       & & MSE & 0.0131 & 0.2358 & 0.0119 & 0.1583 & 0.0182 & 0.2844 \\
     & 3 & Avg & 0.5231 & 1.5892 & 0.5213 & 1.5658 & 0.5392 & 1.6428 \\ 
       & & MSE & 0.0118 & 0.2316 & 0.0104 & 0.1513 & 0.0156 & 0.2955 \\ 
   \hline
\end{tabular}
\end{table}
\begin{table}[ht] 
\centering
\caption{Average values and the corresponding MSEs of the estimators NR, EM and SEM when $n=50$.}
\label{NR:n50}
\begin{tabular}{ccccccccc}
\multicolumn{9}{c}{$m=25$}\\
  \hline
 &&& \multicolumn{2}{c}{NR} & \multicolumn{2}{c}{EM} & \multicolumn{2}{c}{SEM} \\ 
 \cmidrule{4-9} 
 $T$& $R$&&$\aa$&$\bb$&$\aa$&$\bb$&$\aa$&$\bb$ \\ 
 \hline
     & 1 & Avg & 0.5300 & 1.6927 & 0.5137 & 1.5562 & 0.5362 & 1.6928 \\ 
       & & MSE & 0.0127 & 0.3926 & 0.0055 & 0.0961 & 0.0129 & 0.3602 \\ 
0.21 & 2 & Avg & 0.5225 & 1.6671 & 0.5137 & 1.5427 & 0.5462 & 1.7300 \\ 
       & & MSE & 0.0155 & 0.4503 & 0.0074 & 0.0776 & 0.0189 & 0.4644 \\ 
     & 3 & Avg & 0.5233 & 1.6345 & 0.5147 & 1.5464 & 0.5372 & 1.6791 \\ 
       & & MSE & 0.0119 & 0.3016 & 0.0056 & 0.0862 & 0.0142 & 0.3311 \\ 
       \hline 
     & 1 & Avg & 0.5370 & 1.7310 & 0.5155 & 1.5745 & 0.5402 & 1.7216 \\
       & & MSE & 0.0125 & 0.3736 & 0.0056 & 0.0953 & 0.0135 & 0.3737 \\  
1.15 & 2 & Avg & 0.5165 & 1.5796 & 0.5137 & 1.5473 & 0.5297 & 1.6098 \\ 
       & & MSE & 0.0086 & 0.1633 & 0.0074 & 0.0787 & 0.0119 & 0.1964 \\ 
     & 3 & Avg & 0.5328 & 1.6759 & 0.5199 & 1.5805 & 0.5382 & 1.6740 \\  
       & & MSE & 0.0103 & 0.2783 & 0.0066 & 0.1109 & 0.0127 & 0.2895 \\ 
   \hline
\end{tabular}
\begin{tabular}{ccccccccc}
\multicolumn{9}{c}{$m=40$}\\
  \hline
     & 1 & Avg & 0.5180 & 1.5921 & 0.5122 & 1.5399 & 0.5304 & 1.6373 \\
       & & MSE & 0.0104 & 0.2087 & 0.0053 & 0.0832 & 0.0116 & 0.2373 \\ 
0.21 & 2 & Avg & 0.5193 & 1.6038 & 0.5145 & 1.5449 & 0.5337 & 1.6530 \\ 
       & & MSE & 0.0115 & 0.2607 & 0.0060 & 0.0902 & 0.0135 & 0.3019 \\ 
     & 3 & Avg & 0.5214 & 1.6289 & 0.5140 & 1.5511 & 0.5367 & 1.6830 \\ 
       & & MSE & 0.0118 & 0.2818 & 0.0057 & 0.0893 & 0.0138 & 0.3341 \\ 
       \hline 
     & 1 & Avg & 0.5151 & 1.5606 & 0.5137 & 1.5537 & 0.5210 & 1.5771 \\ 
       & & MSE & 0.0058 & 0.0971 & 0.0055 & 0.0785 & 0.0070 & 0.1141 \\
1.15 & 2 & Avg & 0.5140 & 1.5469 & 0.5134 & 1.5393 & 0.5219 & 1.5654 \\ 
       & & MSE & 0.0060 & 0.0879 & 0.0057 & 0.0642 & 0.0077 & 0.1071 \\
     & 3 & Avg & 0.5154 & 1.5536 & 0.5147 & 1.5457 & 0.5224 & 1.5704 \\ 
       & & MSE & 0.0061 & 0.0922 & 0.0057 & 0.0674 & 0.0077 & 0.1109 \\ 
   \hline
\end{tabular}
\end{table}
\begin{sidewaystable} 
  \centering
\caption{Average values and the corresponding MSEs of the Bayes estimators with T-K approximation when $n=30$ and $m=15$.}
\label{TK:n30m15}
\begin{tabular}{ccccccccccccc}
\multicolumn{13}{c}{Informative Priors}\\
  \hline
 &&& \multicolumn{2}{c}{SEL} & \multicolumn{2}{c}{LINEX} & \multicolumn{2}{c}{LINEX} & \multicolumn{2}{c}{GEL} &\multicolumn{2}{c}{GEL}\\
 &&& & & \multicolumn{2}{c}{$\nu=-0.5$} & \multicolumn{2}{c}{$\nu=0.5$} & \multicolumn{2}{c}{$\kappa=-0.5$} &\multicolumn{2}{c}{$\kappa=0.5$}\\
 \cmidrule{4-13} 
 $T$& $R$&&$\aa$&$\bb$&$\aa$&$\bb$&$\aa$&$\bb$&$\aa$&$\bb$&$\aa$&$\bb$ \\ 
 \hline
     & 1 & Avg & 0.5208 & 1.5714 & 0.5217 & 1.5890 & 0.5198 & 1.5537 & 0.5189 & 1.5603 & 0.5152 & 1.5380 \\  
       & & MSE & 0.0017 & 0.0247 & 0.0002 & 0.0034 & 0.0002 & 0.0028 & 0.0007 & 0.0012 & 0.0007 & 0.0011 \\ 
0.21 & 2 & Avg & 0.5207 & 1.5726 & 0.5218 & 1.5928 & 0.5196 & 1.5526 & 0.5187 & 1.5599 & 0.5146 & 1.5345 \\ 
       & & MSE & 0.0016 & 0.0200 & 0.0002 & 0.0029 & 0.0002 & 0.0022 & 0.0007 & 0.0009 & 0.0007 & 0.0008 \\  
     & 3 & Avg & 0.5205 & 1.5706 & 0.5215 & 1.5895 & 0.5195 & 1.5518 & 0.5186 & 1.5587 & 0.5148 & 1.5348 \\ 
       & & MSE & 0.0017 & 0.0215 & 0.0002 & 0.0030 & 0.0002 & 0.0024 & 0.0007 & 0.0010 & 0.0007 & 0.0009 \\
       \hline 
     & 1 & Avg & 0.5205 & 1.5776 & 0.5214 & 1.5949 & 0.5195 & 1.5603 & 0.5186 & 1.5668 & 0.5149 & 1.5452 \\ 
       & & MSE & 0.0018 & 0.0256 & 0.0002 & 0.0036 & 0.0002 & 0.0029 & 0.0008 & 0.0012 & 0.0007 & 0.0011 \\  
1.15 & 2 & Avg & 0.5188 & 1.5640 & 0.5198 & 1.5813 & 0.5177 & 1.5469 & 0.5168 & 1.5531 & 0.5129 & 1.5313 \\  
       & & MSE & 0.0018 & 0.0261 & 0.0002 & 0.0036 & 0.0002 & 0.0030 & 0.0008 & 0.0013 & 0.0007 & 0.0012 \\ 
     & 3 & Avg & 0.5215 & 1.5691 & 0.5225 & 1.5857 & 0.5205 & 1.5526 & 0.5197 & 1.5587 & 0.5160 & 1.5379 \\ 
       & & MSE & 0.0019 & 0.0284 & 0.0002 & 0.0039 & 0.0002 & 0.0032 & 0.0008 & 0.0014 & 0.0008 & 0.0013 \\ 
   \hline
\end{tabular}
\begin{tabular}{ccccccccccccc}
\multicolumn{13}{c}{Non-Informative Priors}\\
  \hline
   &&& \multicolumn{2}{c}{SEL} & \multicolumn{2}{c}{LINEX} & \multicolumn{2}{c}{LINEX} & \multicolumn{2}{c}{GEL} &\multicolumn{2}{c}{GEL}\\
 &&& & & \multicolumn{2}{c}{$\nu=-0.5$} & \multicolumn{2}{c}{$\nu=0.5$} & \multicolumn{2}{c}{$\kappa=-0.5$} &\multicolumn{2}{c}{$\kappa=0.5$}\\
 \cmidrule{4-13} 
  $T$& $R$&&$\aa$&$\bb$&$\aa$&$\bb$&$\aa$&$\bb$&$\aa$&$\bb$&$\aa$&$\bb$ \\ 
 \hline
     & 1 & Avg & 0.5399 & 1.8288 & 0.5320 & 1.8595 & 0.5226 & 1.5580 & 0.5425 & 1.9302 & 0.5248 & 1.7399 \\ 
       & & MSE & 0.0206 & 0.7570 & 0.0022 & 0.0787 & 0.0021 & 0.0364 & 0.0087 & 0.0277 & 0.0090 & 0.0285 \\
0.21 & 2 & Avg & 0.5417 & 1.8412 & 0.5246 & 1.8695 & 0.5116 & 1.4614 & 0.5352 & 1.7773 & 0.5108 & 1.5559 \\ 
       & & MSE & 0.0258 & 0.9058 & 0.0025 & 0.0972 & 0.0023 & 0.0436 & 0.0109 & 0.0372 & 0.0112 & 0.0367 \\  
     & 3 & Avg & 0.5289 & 1.7707 & 0.5246 & 1.8620 & 0.5144 & 1.5292 & 0.5254 & 1.7517 & 0.5063 & 1.5805 \\ 
       & & MSE & 0.0200 & 0.6477 & 0.0022 & 0.0837 & 0.0021 & 0.0439 & 0.0089 & 0.0255 & 0.0092 & 0.0264 \\
       \hline 
     & 1 & Avg & 0.5505 & 1.8844 & 0.5426 & 1.9084 & 0.5339 & 1.6149 & 0.5532 & 1.9776 & 0.5370 & 1.7977 \\ 
       & & MSE & 0.0204 & 0.7696 & 0.0022 & 0.0753 & 0.0021 & 0.0332 & 0.0081 & 0.0236 & 0.0083 & 0.0241 \\
1.15 & 2 & Avg & 0.5263 & 1.6340 & 0.5286 & 1.6998 & 0.5216 & 1.5586 & 0.5207 & 1.6045 & 0.5079 & 1.5266 \\  
       & & MSE & 0.0153 & 0.3325 & 0.0019 & 0.0422 & 0.0019 & 0.0359 & 0.0064 & 0.0139 & 0.0065 & 0.0144 \\ 
     & 3 & Avg & 0.5451 & 1.7624 & 0.5437 & 1.8058 & 0.5370 & 1.6225 & 0.5430 & 1.7819 & 0.5308 & 1.6835 \\ 
       & & MSE & 0.0173 & 0.5539 & 0.0020 & 0.0585 & 0.0020 & 0.0403 & 0.0066 & 0.0177 & 0.0068 & 0.0187 \\ 
   \hline
\end{tabular}
\end{sidewaystable}
\begin{sidewaystable} 
  \centering
\caption{Average values and the corresponding MSEs of the Bayes estimators with T-K approximation when $n=30$ and $m=20$.}
\label{TK:n30m20}
\begin{tabular}{ccccccccccccc}
\multicolumn{13}{c}{Informative Priors}\\
  \hline
 &&& \multicolumn{2}{c}{SEL} & \multicolumn{2}{c}{LINEX} & \multicolumn{2}{c}{LINEX} & \multicolumn{2}{c}{GEL} &\multicolumn{2}{c}{GEL}\\
 &&& & & \multicolumn{2}{c}{$\nu=-0.5$} & \multicolumn{2}{c}{$\nu=0.5$} & \multicolumn{2}{c}{$\kappa=-0.5$} &\multicolumn{2}{c}{$\kappa=0.5$}\\
 \cmidrule{4-13} 
 $T$& $R$&&$\aa$&$\bb$&$\aa$&$\bb$&$\aa$&$\bb$&$\aa$&$\bb$&$\aa$&$\bb$ \\ 
 \hline
     & 1 & Avg & 0.5219 & 1.5643 & 0.5229 & 1.5812 & 0.5209 & 1.5473 & 0.5201 & 1.5535 & 0.5163 & 1.5320 \\  
       & & MSE & 0.0018 & 0.0236 & 0.0002 & 0.0032 & 0.0002 & 0.0027 & 0.0008 & 0.0011 & 0.0007 & 0.0011 \\ 
0.21 & 2 & Avg & 0.5198 & 1.5742 & 0.5209 & 1.5932 & 0.5188 & 1.5553 & 0.5179 & 1.5623 & 0.5139 & 1.5383 \\ 
       & & MSE & 0.0017 & 0.0232 & 0.0002 & 0.0032 & 0.0002 & 0.0026 & 0.0007 & 0.0011 & 0.0007 & 0.0010 \\  
     & 3 & Avg & 0.5201 & 1.5712 & 0.5211 & 1.5901 & 0.5191 & 1.5524 & 0.5182 & 1.5593 & 0.5144 & 1.5354 \\ 
       & & MSE & 0.0017 & 0.0218 & 0.0002 & 0.0031 & 0.0002 & 0.0024 & 0.0007 & 0.0010 & 0.0007 & 0.0009 \\
       \hline 
     & 1 & Avg & 0.5213 & 1.5692 & 0.5222 & 1.5840 & 0.5203 & 1.5544 & 0.5195 & 1.5599 & 0.5160 & 1.5413 \\ 
       & & MSE & 0.0018 & 0.0280 & 0.0002 & 0.0038 & 0.0002 & 0.0032 & 0.0008 & 0.0013 & 0.0008 & 0.0013 \\  
1.15 & 2 & Avg & 0.5206 & 1.5584 & 0.5215 & 1.5738 & 0.5196 & 1.5430 & 0.5187 & 1.5486 & 0.5150 & 1.5290 \\  
       & & MSE & 0.0019 & 0.0280 & 0.0002 & 0.0037 & 0.0002 & 0.0033 & 0.0008 & 0.0014 & 0.0008 & 0.0013 \\ 
     & 3 & Avg & 0.5193 & 1.5513 & 0.5202 & 1.5667 & 0.5183 & 1.5361 & 0.5175 & 1.5415 & 0.5139 & 1.5219 \\ 
       & & MSE & 0.0019 & 0.0274 & 0.0002 & 0.0036 & 0.0002 & 0.0032 & 0.0008 & 0.0014 & 0.0008 & 0.0013 \\ 
   \hline
\end{tabular}
\begin{tabular}{ccccccccccccc}
\multicolumn{13}{c}{Non-Informative Priors}\\
  \hline
   &&& \multicolumn{2}{c}{SEL} & \multicolumn{2}{c}{LINEX} & \multicolumn{2}{c}{LINEX} & \multicolumn{2}{c}{GEL} &\multicolumn{2}{c}{GEL}\\
 &&& & & \multicolumn{2}{c}{$\nu=-0.5$} & \multicolumn{2}{c}{$\nu=0.5$} & \multicolumn{2}{c}{$\kappa=-0.5$} &\multicolumn{2}{c}{$\kappa=0.5$}\\
 \cmidrule{4-13} 
  $T$& $R$&&$\aa$&$\bb$&$\aa$&$\bb$&$\aa$&$\bb$&$\aa$&$\bb$&$\aa$&$\bb$ \\ 
 \hline
     & 1 & Avg & 0.5357 & 1.7173 & 0.5334 & 1.7833 & 0.5244 & 1.5617 & 0.5299 & 1.6836 & 0.5137 & 1.5703 \\ 
       & & MSE & 0.0202 & 0.4813 & 0.0022 & 0.0594 & 0.0022 & 0.0389 & 0.0078 & 0.0178 & 0.0080 & 0.0181 \\
0.21 & 2 & Avg & 0.5316 & 1.7861 & 0.5238 & 1.8553 & 0.5130 & 1.5358 & 0.5265 & 1.7418 & 0.5064 & 1.5767 \\ 
       & & MSE & 0.0222 & 0.6907 & 0.0023 & 0.0824 & 0.0022 & 0.0457 & 0.0094 & 0.0258 & 0.0097 & 0.0263 \\  
     & 3 & Avg & 0.5295 & 1.7864 & 0.5216 & 1.8587 & 0.5116 & 1.5283 & 0.5233 & 1.7309 & 0.5043 & 1.5660 \\ 
       & & MSE & 0.0212 & 0.7066 & 0.0022 & 0.0846 & 0.0021 & 0.0451 & 0.0091 & 0.0259 & 0.0093 & 0.0260 \\
       \hline 
     & 1 & Avg & 0.5425 & 1.7180 & 0.5432 & 1.7634 & 0.5370 & 1.6293 & 0.5383 & 1.7009 & 0.5273 & 1.6338 \\ 
       & & MSE & 0.0151 & 0.3740 & 0.0018 & 0.0428 & 0.0018 & 0.0330 & 0.0055 & 0.0112 & 0.0056 & 0.0115 \\
1.15 & 2 & Avg & 0.5283 & 1.5985 & 0.5309 & 1.6459 & 0.5250 & 1.5488 & 0.5231 & 1.5715 & 0.5123 & 1.5158 \\  
       & & MSE & 0.0131 & 0.2337 & 0.0016 & 0.0295 & 0.0016 & 0.0276 & 0.0053 & 0.0101 & 0.0054 & 0.0103 \\ 
     & 3 & Avg & 0.5211 & 1.5766 & 0.5235 & 1.6232 & 0.5181 & 1.5257 & 0.5163 & 1.5521 & 0.5062 & 1.4954 \\ 
       & & MSE & 0.0117 & 0.2284 & 0.0015 & 0.0288 & 0.0014 & 0.0262 & 0.0051 & 0.0107 & 0.0052 & 0.0109 \\ 
   \hline
\end{tabular}
\end{sidewaystable}
\begin{sidewaystable} 
  \centering
\caption{Average values and the corresponding MSEs of the Bayes estimators with MCMC method when $n=30$ and $m=15$.}
\label{MCMC:n30m15}
\begin{tabular}{ccccccccccccc}
\multicolumn{13}{c}{Informative Priors}\\
  \hline
 &&& \multicolumn{2}{c}{SEL} & \multicolumn{2}{c}{LINEX} & \multicolumn{2}{c}{LINEX} & \multicolumn{2}{c}{GEL} &\multicolumn{2}{c}{GEL}\\
 &&& & & \multicolumn{2}{c}{$\nu=-0.5$} & \multicolumn{2}{c}{$\nu=0.5$} & \multicolumn{2}{c}{$\kappa=-0.5$} &\multicolumn{2}{c}{$\kappa=0.5$}\\
 \cmidrule{4-13} 
 $T$& $R$&&$\aa$&$\bb$&$\aa$&$\bb$&$\aa$&$\bb$&$\aa$&$\bb$&$\aa$&$\bb$ \\ 
 \hline
     & 1 & Avg & 0.5208 & 1.5717 & 0.5217 & 1.5898 & 0.5198 & 1.5543 & 0.5189 & 1.5606 & 0.5152 & 1.5385 \\  
       & & MSE & 0.0017 & 0.0262 & 0.0002 & 0.0042 & 0.0002 & 0.0030 & 0.0008 & 0.0015 & 0.0008 & 0.0012 \\ 
0.21 & 2 & Avg & 0.5206 & 1.5725 & 0.5216 & 1.5928 & 0.5195 & 1.5528 & 0.5185 & 1.5599 & 0.5144 & 1.5345 \\ 
       & & MSE & 0.0016 & 0.0204 & 0.0002 & 0.0032 & 0.0002 & 0.0023 & 0.0008 & 0.0012 & 0.0007 & 0.0009 \\  
     & 3 & Avg & 0.5205 & 1.5710 & 0.5215 & 1.5902 & 0.5195 & 1.5524 & 0.5186 & 1.5591 & 0.5147 & 1.5352 \\ 
       & & MSE & 0.0017 & 0.0225 & 0.0002 & 0.0036 & 0.0002 & 0.0025 & 0.0008 & 0.0013 & 0.0008 & 0.0010 \\
       \hline 
     & 1 & Avg & 0.5204 & 1.5778 & 0.5214 & 1.5953 & 0.5195 & 1.5609 & 0.5186 & 1.5671 & 0.5149 & 1.5456 \\ 
       & & MSE & 0.0018 & 0.0263 & 0.0002 & 0.0039 & 0.0002 & 0.0030 & 0.0008 & 0.0015 & 0.0008 & 0.0012 \\  
1.15 & 2 & Avg & 0.5187 & 1.5639 & 0.5197 & 1.5812 & 0.5177 & 1.5470 & 0.5168 & 1.5530 & 0.5129 & 1.5314 \\  
       & & MSE & 0.0018 & 0.0263 & 0.0002 & 0.0038 & 0.0002 & 0.0030 & 0.0008 & 0.0015 & 0.0008 & 0.0012 \\ 
     & 3 & Avg & 0.5215 & 1.5690 & 0.5225 & 1.5857 & 0.5206 & 1.5528 & 0.5197 & 1.5587 & 0.5160 & 1.5379 \\ 
       & & MSE & 0.0019 & 0.0287 & 0.0002 & 0.0042 & 0.0002 & 0.0033 & 0.0009 & 0.0016 & 0.0009 & 0.0014 \\ 
   \hline
\end{tabular}
\begin{tabular}{ccccccccccccc}
\multicolumn{13}{c}{Non-Informative Priors}\\
  \hline
   &&& \multicolumn{2}{c}{SEL} & \multicolumn{2}{c}{LINEX} & \multicolumn{2}{c}{LINEX} & \multicolumn{2}{c}{GEL} &\multicolumn{2}{c}{GEL}\\
 &&& & & \multicolumn{2}{c}{$\nu=-0.5$} & \multicolumn{2}{c}{$\nu=0.5$} & \multicolumn{2}{c}{$\kappa=-0.5$} &\multicolumn{2}{c}{$\kappa=0.5$}\\
 \cmidrule{4-13} 
  $T$& $R$&&$\aa$&$\bb$&$\aa$&$\bb$&$\aa$&$\bb$&$\aa$&$\bb$&$\aa$&$\bb$ \\ 
 \hline
     & 1 & Avg & 0.5273 & 1.7244 & 0.5320 & 1.9407 & 0.5227 & 1.5963 & 0.5190 & 1.6569 & 0.5020 & 1.5310 \\ 
       & & MSE & 0.0163 & 0.4777 & 0.0022 & 0.1824 & 0.0020 & 0.0434 & 0.0073 & 0.0296 & 0.0069 & 0.0152 \\
0.21 & 2 & Avg & 0.5170 & 1.7150 & 0.5235 & 2.0844 & 0.5107 & 1.5342 & 0.5054 & 1.6121 & 0.4818 & 1.4226 \\ 
       & & MSE & 0.0174 & 0.5795 & 0.0024 & 0.2607 & 0.0021 & 0.0512 & 0.0087 & 0.0438 & 0.0083 & 0.0242 \\  
     & 3 & Avg & 0.5178 & 1.7089 & 0.5229 & 1.9536 & 0.5129 & 1.5662 & 0.5086 & 1.6291 & 0.4899 & 1.4801 \\ 
       & & MSE & 0.0163 & 0.4933 & 0.0022 & 0.1774 & 0.0020 & 0.0469 & 0.0077 & 0.0313 & 0.0074 & 0.0178 \\
       \hline 
     & 1 & Avg & 0.5410 & 1.7834 & 0.5454 & 1.9879 & 0.5367 & 1.6582 & 0.5334 & 1.7202 & 0.5179 & 1.6014 \\ 
       & & MSE & 0.0176 & 0.4799 & 0.0024 & 0.1840 & 0.0022 & 0.0431 & 0.0073 & 0.0278 & 0.0068 & 0.0124 \\
1.15 & 2 & Avg & 0.5247 & 1.6263 & 0.5282 & 1.7012 & 0.5213 & 1.5614 & 0.5184 & 1.5884 & 0.5056 & 1.5124 \\  
       & & MSE & 0.0150 & 0.3141 & 0.0020 & 0.0640 & 0.0019 & 0.0365 & 0.0065 & 0.0155 & 0.0063 & 0.0124 \\ 
     & 3 & Avg & 0.5412 & 1.7179 & 0.5446 & 1.8252 & 0.5379 & 1.6365 & 0.5353 & 1.6760 & 0.5233 & 1.5936 \\ 
       & & MSE & 0.0162 & 0.4243 & 0.0022 & 0.1179 & 0.0020 & 0.0453 & 0.0067 & 0.0205 & 0.0064 & 0.0129 \\ 
   \hline
\end{tabular}
\end{sidewaystable}
\begin{sidewaystable} 
  \centering
\caption{Average values and the corresponding MSEs of the Bayes estimators with MCMC method when $n=30$ and $m=20$.}
\label{MCMC:n30m20}
\begin{tabular}{ccccccccccccc}
\multicolumn{13}{c}{Informative Priors}\\
  \hline
 &&& \multicolumn{2}{c}{SEL} & \multicolumn{2}{c}{LINEX} & \multicolumn{2}{c}{LINEX} & \multicolumn{2}{c}{GEL} &\multicolumn{2}{c}{GEL}\\
 &&& & & \multicolumn{2}{c}{$\nu=-0.5$} & \multicolumn{2}{c}{$\nu=0.5$} & \multicolumn{2}{c}{$\kappa=-0.5$} &\multicolumn{2}{c}{$\kappa=0.5$}\\
 \cmidrule{4-13} 
 $T$& $R$&&$\aa$&$\bb$&$\aa$&$\bb$&$\aa$&$\bb$&$\aa$&$\bb$&$\aa$&$\bb$ \\ 
 \hline
     & 1 & Avg & 0.5218 & 1.5638 & 0.5228 & 1.5809 & 0.5208 & 1.5472 & 0.5200 & 1.5531 & 0.5162 & 1.5317 \\  
       & & MSE & 0.0018 & 0.0239 & 0.0002 & 0.0035 & 0.0002 & 0.0028 & 0.0009 & 0.0013 & 0.0008 & 0.0011 \\ 
0.21 & 2 & Avg & 0.5199 & 1.5744 & 0.5209 & 1.5935 & 0.5188 & 1.5559 & 0.5179 & 1.5625 & 0.5139 & 1.5387 \\ 
       & & MSE & 0.0017 & 0.0237 & 0.0002 & 0.0036 & 0.0002 & 0.0027 & 0.0008 & 0.0014 & 0.0008 & 0.0011 \\  
     & 3 & Avg & 0.5201 & 1.5707 & 0.5211 & 1.5898 & 0.5191 & 1.5523 & 0.5182 & 1.5589 & 0.5143 & 1.5352 \\ 
       & & MSE & 0.0017 & 0.0220 & 0.0002 & 0.0033 & 0.0002 & 0.0025 & 0.0008 & 0.0013 & 0.0008 & 0.0010 \\
       \hline 
     & 1 & Avg & 0.5213 & 1.5690 & 0.5222 & 1.5839 & 0.5203 & 1.5544 & 0.5195 & 1.5597 & 0.5160 & 1.5412 \\ 
       & & MSE & 0.0019 & 0.0284 & 0.0002 & 0.0041 & 0.0002 & 0.0033 & 0.0009 & 0.0016 & 0.0008 & 0.0013 \\  
1.15 & 2 & Avg & 0.5206 & 1.5582 & 0.5216 & 1.5737 & 0.5196 & 1.5431 & 0.5188 & 1.5485 & 0.5151 & 1.5290 \\  
       & & MSE & 0.0019 & 0.0282 & 0.0002 & 0.0040 & 0.0002 & 0.0033 & 0.0009 & 0.0016 & 0.0009 & 0.0014 \\ 
     & 3 & Avg & 0.5192 & 1.5514 & 0.5201 & 1.5669 & 0.5183 & 1.5364 & 0.5174 & 1.5417 & 0.5139 & 1.5221 \\ 
       & & MSE & 0.0019 & 0.0277 & 0.0002 & 0.0039 & 0.0002 & 0.0033 & 0.0009 & 0.0015 & 0.0008 & 0.0014 \\ 
   \hline
\end{tabular}
\begin{tabular}{ccccccccccccc}
\multicolumn{13}{c}{Non-Informative Priors}\\
  \hline
   &&& \multicolumn{2}{c}{SEL} & \multicolumn{2}{c}{LINEX} & \multicolumn{2}{c}{LINEX} & \multicolumn{2}{c}{GEL} &\multicolumn{2}{c}{GEL}\\
 &&& & & \multicolumn{2}{c}{$\nu=-0.5$} & \multicolumn{2}{c}{$\nu=0.5$} & \multicolumn{2}{c}{$\kappa=-0.5$} &\multicolumn{2}{c}{$\kappa=0.5$}\\
 \cmidrule{4-13} 
  $T$& $R$&&$\aa$&$\bb$&$\aa$&$\bb$&$\aa$&$\bb$&$\aa$&$\bb$&$\aa$&$\bb$ \\ 
 \hline
     & 1 & Avg & 0.5303 & 1.6921 & 0.5348 & 1.8457 & 0.5259 & 1.5907 & 0.5224 & 1.6362 & 0.5064 & 1.5308 \\ 
       & & MSE & 0.0178 & 0.4075 & 0.0024 & 0.1214 & 0.0022 & 0.0421 & 0.0076 & 0.0228 & 0.0072 & 0.0142 \\
0.21 & 2 & Avg & 0.5182 & 1.7196 & 0.5236 & 1.9673 & 0.5129 & 1.5778 & 0.5085 & 1.6406 & 0.4889 & 1.4930 \\ 
       & & MSE & 0.0169 & 0.5114 & 0.0023 & 0.1816 & 0.0020 & 0.0493 & 0.0080 & 0.0322 & 0.0076 & 0.0184 \\  
     & 3 & Avg & 0.5162 & 1.7124 & 0.5213 & 1.9567 & 0.5113 & 1.5689 & 0.5071 & 1.6326 & 0.4885 & 1.4838 \\ 
       & & MSE & 0.0163 & 0.5117 & 0.0022 & 0.1811 & 0.0020 & 0.0487 & 0.0078 & 0.0323 & 0.0075 & 0.0186 \\
       \hline 
     & 1 & Avg & 0.5406 & 1.7004 & 0.5436 & 1.7762 & 0.5375 & 1.6393 & 0.5352 & 1.6684 & 0.5243 & 1.6058 \\ 
       & & MSE & 0.0145 & 0.3223 & 0.0019 & 0.0786 & 0.0018 & 0.0361 & 0.0059 & 0.0149 & 0.0056 & 0.0098 \\
1.15 & 2 & Avg & 0.5277 & 1.5945 & 0.5307 & 1.6447 & 0.5248 & 1.5490 & 0.5224 & 1.5671 & 0.5116 & 1.5121 \\  
       & & MSE & 0.0130 & 0.2263 & 0.0017 & 0.0407 & 0.0016 & 0.0275 & 0.0056 & 0.0111 & 0.0054 & 0.0095 \\ 
     & 3 & Avg & 0.5206 & 1.5734 & 0.5232 & 1.6243 & 0.5179 & 1.5273 & 0.5156 & 1.5454 & 0.5055 & 1.4893 \\ 
       & & MSE & 0.0116 & 0.2250 & 0.0015 & 0.0404 & 0.0014 & 0.0271 & 0.0052 & 0.0114 & 0.0051 & 0.0099 \\ 
   \hline
\end{tabular}
\end{sidewaystable}

\begin{sidewaystable} 
  \centering
\caption{Average values and the corresponding MSEs of the Bayes estimators with T-K approximation when $n=50$ and $m=25$.}
\label{TK:n50m25}
\begin{tabular}{ccccccccccccc}
\multicolumn{13}{c}{Informative Priors}\\
  \hline
 &&& \multicolumn{2}{c}{SEL} & \multicolumn{2}{c}{LINEX} & \multicolumn{2}{c}{LINEX} & \multicolumn{2}{c}{GEL} &\multicolumn{2}{c}{GEL}\\
 &&& & & \multicolumn{2}{c}{$\nu=-0.5$} & \multicolumn{2}{c}{$\nu=0.5$} & \multicolumn{2}{c}{$\kappa=-0.5$} &\multicolumn{2}{c}{$\kappa=0.5$}\\
 \cmidrule{4-13} 
 $T$& $R$&&$\aa$&$\bb$&$\aa$&$\bb$&$\aa$&$\bb$&$\aa$&$\bb$&$\aa$&$\bb$ \\ 
 \hline
     & 1 & Avg & 0.5189 & 1.5698 & 0.5196 & 1.5846 & 0.5180 & 1.5547 & 0.5173 & 1.5603 & 0.5142 & 1.5415 \\  
       & & MSE & 0.0017 & 0.0267 & 0.0002 & 0.0036 & 0.0002 & 0.0031 & 0.0008 & 0.0013 & 0.0007 & 0.0012 \\ 
0.21 & 2 & Avg & 0.5196 & 1.5682 & 0.5205 & 1.5861 & 0.5186 & 1.5503 & 0.5178 & 1.5569 & 0.5142 & 1.5342 \\ 
       & & MSE & 0.0017 & 0.0230 & 0.0002 & 0.0032 & 0.0002 & 0.0026 & 0.0007 & 0.0011 & 0.0007 & 0.0010 \\  
     & 3 & Avg & 0.5193 & 1.5651 & 0.5201 & 1.5806 & 0.5184 & 1.5495 & 0.5177 & 1.5553 & 0.5145 & 1.5355 \\ 
       & & MSE & 0.0017 & 0.0249 & 0.0002 & 0.0033 & 0.0002 & 0.0029 & 0.0008 & 0.0012 & 0.0007 & 0.0011 \\
       \hline 
     & 1 & Avg & 0.5191 & 1.5771 & 0.5199 & 1.5917 & 0.5183 & 1.5624 & 0.5176 & 1.5680 & 0.5146 & 1.5497 \\ 
       & & MSE & 0.0018 & 0.0273 & 0.0002 & 0.0037 & 0.0002 & 0.0031 & 0.0008 & 0.0013 & 0.0007 & 0.0012 \\  
1.15 & 2 & Avg & 0.5177 & 1.5557 & 0.5185 & 1.5698 & 0.5168 & 1.5417 & 0.5160 & 1.5468 & 0.5128 & 1.5288 \\  
       & & MSE & 0.0018 & 0.0278 & 0.0002 & 0.0037 & 0.0002 & 0.0033 & 0.0008 & 0.0014 & 0.0008 & 0.0013 \\ 
     & 3 & Avg & 0.5197 & 1.5673 & 0.5204 & 1.5810 & 0.5189 & 1.5534 & 0.5182 & 1.5586 & 0.5151 & 1.5412 \\ 
       & & MSE & 0.0018 & 0.0284 & 0.0002 & 0.0038 & 0.0002 & 0.0033 & 0.0008 & 0.0014 & 0.0008 & 0.0013 \\ 
   \hline
\end{tabular}
\begin{tabular}{ccccccccccccc}
\multicolumn{13}{c}{Non-Informative Priors}\\
  \hline
   &&& \multicolumn{2}{c}{SEL} & \multicolumn{2}{c}{LINEX} & \multicolumn{2}{c}{LINEX} & \multicolumn{2}{c}{GEL} &\multicolumn{2}{c}{GEL}\\
 &&& & & \multicolumn{2}{c}{$\nu=-0.5$} & \multicolumn{2}{c}{$\nu=0.5$} & \multicolumn{2}{c}{$\kappa=-0.5$} &\multicolumn{2}{c}{$\kappa=0.5$}\\
 \cmidrule{4-13} 
  $T$& $R$&&$\aa$&$\bb$&$\aa$&$\bb$&$\aa$&$\bb$&$\aa$&$\bb$&$\aa$&$\bb$ \\ 
 \hline
     & 1 & Avg & 0.5288 & 1.7240 & 0.5283 & 1.7734 & 0.5230 & 1.6094 & 0.5256 & 1.7062 & 0.5157 & 1.6251 \\ 
       & & MSE & 0.0125 & 0.4397 & 0.0015 & 0.0509 & 0.0015 & 0.0363 & 0.0052 & 0.0142 & 0.0053 & 0.0144 \\
0.21 & 2 & Avg & 0.5248 & 1.7058 & 0.5234 & 1.7885 & 0.5158 & 1.5500 & 0.5190 & 1.6571 & 0.5048 & 1.5365 \\ 
       & & MSE & 0.0156 & 0.4929 & 0.0018 & 0.0645 & 0.0018 & 0.0433 & 0.0067 & 0.0198 & 0.0068 & 0.0201 \\  
     & 3 & Avg & 0.5223 & 1.6580 & 0.5235 & 1.7186 & 0.5180 & 1.5688 & 0.5175 & 1.6224 & 0.5072 & 1.5450 \\ 
       & & MSE & 0.0117 & 0.3217 & 0.0014 & 0.0421 & 0.0014 & 0.0330 & 0.0050 & 0.0129 & 0.0051 & 0.0129 \\
       \hline 
     & 1 & Avg & 0.5358 & 1.7599 & 0.5360 & 1.8176 & 0.5309 & 1.6538 & 0.5326 & 1.7421 & 0.5233 & 1.6646 \\ 
       & & MSE & 0.0123 & 0.4209 & 0.0015 & 0.0512 & 0.0015 & 0.0360 & 0.0048 & 0.0122 & 0.0049 & 0.0123 \\
1.15 & 2 & Avg & 0.5178 & 1.5715 & 0.5197 & 1.6080 & 0.5155 & 1.5347 & 0.5139 & 1.5496 & 0.5061 & 1.5058 \\  
       & & MSE & 0.0086 & 0.1622 & 0.0011 & 0.0208 & 0.0011 & 0.0196 & 0.0039 & 0.0076 & 0.0039 & 0.0077 \\ 
     & 3 & Avg & 0.5321 & 1.6862 & 0.5332 & 1.7312 & 0.5288 & 1.6203 & 0.5284 & 1.6617 & 0.5205 & 1.6067 \\ 
       & & MSE & 0.0103 & 0.2985 & 0.0013 & 0.0363 & 0.0012 & 0.0302 & 0.0040 & 0.0092 & 0.0041 & 0.0093 \\ 
   \hline
\end{tabular}
\end{sidewaystable}
\begin{sidewaystable}
  \centering
\caption{Average values and the corresponding MSEs of the Bayes estimators with T-K approximation when $n=50$ and $m=40$.}
\label{TK:n50m40}
\begin{tabular}{ccccccccccccc}
\multicolumn{13}{c}{Informative Priors}\\
  \hline
 &&& \multicolumn{2}{c}{SEL} & \multicolumn{2}{c}{LINEX} & \multicolumn{2}{c}{LINEX} & \multicolumn{2}{c}{GEL} &\multicolumn{2}{c}{GEL}\\
 &&& & & \multicolumn{2}{c}{$\nu=-0.5$} & \multicolumn{2}{c}{$\nu=0.5$} & \multicolumn{2}{c}{$\kappa=-0.5$} &\multicolumn{2}{c}{$\kappa=0.5$}\\
 \cmidrule{4-13} 
 $T$& $R$&&$\aa$&$\bb$&$\aa$&$\bb$&$\aa$&$\bb$&$\aa$&$\bb$&$\aa$&$\bb$ \\ 
 \hline
     & 1 & Avg & 0.5179 & 1.5598 & 0.5186 & 1.5740 & 0.5170 & 1.5455 & 0.5163 & 1.5508 & 0.5132 & 1.5326 \\  
       & & MSE & 0.0017 & 0.0244 & 0.0002 & 0.0032 & 0.0002 & 0.0029 & 0.0007 & 0.0012 & 0.0007 & 0.0011 \\ 
0.21 & 2 & Avg & 0.5188 & 1.5602 & 0.5196 & 1.5756 & 0.5180 & 1.5447 & 0.5172 & 1.5504 & 0.5139 & 1.5308 \\ 
       & & MSE & 0.0017 & 0.0239 & 0.0002 & 0.0032 & 0.0002 & 0.0028 & 0.0007 & 0.0011 & 0.0007 & 0.0011 \\  
     & 3 & Avg & 0.5186 & 1.5667 & 0.5194 & 1.5822 & 0.5178 & 1.5511 & 0.5171 & 1.5569 & 0.5139 & 1.5372 \\ 
       & & MSE & 0.0017 & 0.0249 & 0.0002 & 0.0033 & 0.0002 & 0.0029 & 0.0007 & 0.0012 & 0.0007 & 0.0011 \\
       \hline 
     & 1 & Avg & 0.5155 & 1.5481 & 0.5162 & 1.5581 & 0.5148 & 1.5381 & 0.5142 & 1.5417 & 0.5116 & 1.5289 \\ 
       & & MSE & 0.0017 & 0.0280 & 0.0002 & 0.0036 & 0.0002 & 0.0034 & 0.0008 & 0.0014 & 0.0007 & 0.0014 \\  
1.15 & 2 & Avg & 0.5165 & 1.5471 & 0.5173 & 1.5581 & 0.5158 & 1.5361 & 0.5151 & 1.5400 & 0.5123 & 1.5259 \\ 
       & & MSE & 0.0018 & 0.0263 & 0.0002 & 0.0034 & 0.0002 & 0.0032 & 0.0008 & 0.0013 & 0.0008 & 0.0013 \\ 
     & 3 & Avg & 0.5170 & 1.5499 & 0.5177 & 1.5610 & 0.5163 & 1.5389 & 0.5157 & 1.5428 & 0.5129 & 1.5287 \\ 
       & & MSE & 0.0018 & 0.0268 & 0.0002 & 0.0035 & 0.0002 & 0.0032 & 0.0008 & 0.0013 & 0.0008 & 0.0013 \\ 
   \hline
\end{tabular}
\begin{tabular}{ccccccccccccc}
\multicolumn{13}{c}{Non-Informative Priors}\\
  \hline
   &&& \multicolumn{2}{c}{SEL} & \multicolumn{2}{c}{LINEX} & \multicolumn{2}{c}{LINEX} & \multicolumn{2}{c}{GEL} &\multicolumn{2}{c}{GEL}\\
 &&& & & \multicolumn{2}{c}{$\nu=-0.5$} & \multicolumn{2}{c}{$\nu=0.5$} & \multicolumn{2}{c}{$\kappa=-0.5$} &\multicolumn{2}{c}{$\kappa=0.5$}\\
 \cmidrule{4-13} 
  $T$& $R$&&$\aa$&$\bb$&$\aa$&$\bb$&$\aa$&$\bb$&$\aa$&$\bb$&$\aa$&$\bb$ \\ 
 \hline
     & 1 & Avg & 0.5171 & 1.6121 & 0.5186 & 1.6578 & 0.5137 & 1.5487 & 0.5126 & 1.5821 & 0.5035 & 1.5225 \\  
       & & MSE & 0.0103 & 0.2216 & 0.0013 & 0.0291 & 0.0013 & 0.0247 & 0.0045 & 0.0096 & 0.0046 & 0.0096 \\
0.21 & 2 & Avg & 0.5193 & 1.6293 & 0.5206 & 1.6851 & 0.5149 & 1.5460 & 0.5141 & 1.5915 & 0.5034 & 1.5173 \\ 
       & & MSE & 0.0115 & 0.2803 & 0.0014 & 0.0369 & 0.0014 & 0.0297 & 0.0050 & 0.0119 & 0.0051 & 0.0119 \\  
     & 3 & Avg & 0.5205 & 1.6524 & 0.5219 & 1.7127 & 0.5165 & 1.5683 & 0.5157 & 1.6149 & 0.5055 & 1.5395 \\ 
       & & MSE & 0.0117 & 0.3022 & 0.0014 & 0.0401 & 0.0014 & 0.0324 & 0.0050 & 0.0122 & 0.0051 & 0.0123 \\
       \hline 
     & 1 & Avg & 0.5141 & 1.5577 & 0.5154 & 1.5767 & 0.5128 & 1.5384 & 0.5116 & 1.5462 & 0.5066 & 1.5231 \\ 
       & & MSE & 0.0057 & 0.0976 & 0.0007 & 0.0122 & 0.0007 & 0.0121 & 0.0025 & 0.0043 & 0.0025 & 0.0044 \\
1.15 & 2 & Avg & 0.5135 & 1.5417 & 0.5149 & 1.5628 & 0.5119 & 1.5205 & 0.5107 & 1.5285 & 0.5050 & 1.5021 \\ 
       & & MSE & 0.0060 & 0.0871 & 0.0008 & 0.0111 & 0.0007 & 0.0106 & 0.0027 & 0.0044 & 0.0028 & 0.0044 \\ 
     & 3 & Avg & 0.5146 & 1.5478 & 0.5159 & 1.5693 & 0.5131 & 1.5263 & 0.5119 & 1.5344 & 0.5064 & 1.5077 \\ 
       & & MSE & 0.0061 & 0.0913 & 0.0008 & 0.0116 & 0.0008 & 0.0112 & 0.0027 & 0.0045 & 0.0027 & 0.0045 \\  
   \hline
\end{tabular}
\end{sidewaystable}
\begin{sidewaystable}
  \centering
\caption{Average values and the corresponding MSEs of the Bayes estimators with MCMC Method when $n=50$ and $m=25$.}
\label{MCMC:n50m25}
\begin{tabular}{ccccccccccccc}
\multicolumn{13}{c}{Informative Priors}\\
  \hline
 &&& \multicolumn{2}{c}{SEL} & \multicolumn{2}{c}{LINEX} & \multicolumn{2}{c}{LINEX} & \multicolumn{2}{c}{GEL} &\multicolumn{2}{c}{GEL}\\
 &&& & & \multicolumn{2}{c}{$\nu=-0.5$} & \multicolumn{2}{c}{$\nu=0.5$} & \multicolumn{2}{c}{$\kappa=-0.5$} &\multicolumn{2}{c}{$\kappa=0.5$}\\
 \cmidrule{4-13} 
 $T$& $R$&&$\aa$&$\bb$&$\aa$&$\bb$&$\aa$&$\bb$&$\aa$&$\bb$&$\aa$&$\bb$ \\ 
 \hline
     & 1 & Avg & 0.5188 & 1.5698 & 0.5196 & 1.5849 & 0.5180 & 1.5551 & 0.5173 & 1.5604 & 0.5142 & 1.5416 \\  
       & & MSE & 0.0017 & 0.0271 & 0.0002 & 0.0039 & 0.0002 & 0.0032 & 0.0008 & 0.0015 & 0.0008 & 0.0013 \\ 
0.21 & 2 & Avg & 0.5195 & 1.5683 & 0.5204 & 1.5863 & 0.5186 & 1.5508 & 0.5177 & 1.5570 & 0.5142 & 1.5345 \\ 
       & & MSE & 0.0017 & 0.0235 & 0.0002 & 0.0035 & 0.0002 & 0.0027 & 0.0008 & 0.0013 & 0.0008 & 0.0011 \\  
     & 3 & Avg & 0.5192 & 1.5648 & 0.5200 & 1.5805 & 0.5184 & 1.5495 & 0.5176 & 1.5550 & 0.5145 & 1.5353 \\ 
       & & MSE & 0.0017 & 0.0252 & 0.0002 & 0.0036 & 0.0002 & 0.0029 & 0.0008 & 0.0014 & 0.0008 & 0.0012 \\
       \hline 
     & 1 & Avg & 0.5191 & 1.5769 & 0.5199 & 1.5917 & 0.5183 & 1.5626 & 0.5176 & 1.5679 & 0.5146 & 1.5497 \\ 
       & & MSE & 0.0018 & 0.0276 & 0.0002 & 0.0040 & 0.0002 & 0.0032 & 0.0008 & 0.0015 & 0.0008 & 0.0013 \\  
1.15 & 2 & Avg & 0.5176 & 1.5555 & 0.5185 & 1.5697 & 0.5168 & 1.5416 & 0.5160 & 1.5466 & 0.5128 & 1.5287 \\ 
       & & MSE & 0.0018 & 0.0280 & 0.0002 & 0.0040 & 0.0002 & 0.0033 & 0.0009 & 0.0015 & 0.0008 & 0.0014 \\ 
     & 3 & Avg & 0.5196 & 1.5669 & 0.5204 & 1.5808 & 0.5188 & 1.5533 & 0.5181 & 1.5583 & 0.5151 & 1.5410 \\ 
       & & MSE & 0.0018 & 0.0286 & 0.0002 & 0.0041 & 0.0002 & 0.0033 & 0.0009 & 0.0016 & 0.0008 & 0.0014 \\ 
   \hline
\end{tabular}
\begin{tabular}{ccccccccccccc}
\multicolumn{13}{c}{Non-Informative Priors}\\
  \hline
   &&& \multicolumn{2}{c}{SEL} & \multicolumn{2}{c}{LINEX} & \multicolumn{2}{c}{LINEX} & \multicolumn{2}{c}{GEL} &\multicolumn{2}{c}{GEL}\\
 &&& & & \multicolumn{2}{c}{$\nu=-0.5$} & \multicolumn{2}{c}{$\nu=0.5$} & \multicolumn{2}{c}{$\kappa=-0.5$} &\multicolumn{2}{c}{$\kappa=0.5$}\\
 \cmidrule{4-13} 
  $T$& $R$&&$\aa$&$\bb$&$\aa$&$\bb$&$\aa$&$\bb$&$\aa$&$\bb$&$\aa$&$\bb$ \\ 
 \hline
     & 1 & Avg & 0.5265 & 1.7003 & 0.5291 & 1.7978 & 0.5238 & 1.6267 & 0.5216 & 1.6616 & 0.5119 & 1.5875 \\  
       & & MSE & 0.0119 & 0.3647 & 0.0016 & 0.0929 & 0.0015 & 0.0395 & 0.0052 & 0.0180 & 0.0051 & 0.0117 \\
0.21 & 2 & Avg & 0.5212 & 1.6898 & 0.5251 & 1.8516 & 0.5174 & 1.5805 & 0.5142 & 1.6285 & 0.5001 & 1.5119 \\ 
       & & MSE & 0.0143 & 0.4405 & 0.0019 & 0.1249 & 0.0018 & 0.0459 & 0.0065 & 0.0246 & 0.0063 & 0.0162 \\  
     & 3 & Avg & 0.5219 & 1.6555 & 0.5246 & 1.7439 & 0.5191 & 1.5857 & 0.5168 & 1.6161 & 0.5065 & 1.5402 \\ 
       & & MSE & 0.0117 & 0.3126 & 0.0015 & 0.0702 & 0.0015 & 0.0350 & 0.0052 & 0.0160 & 0.0050 & 0.0117 \\
       \hline 
     & 1 & Avg & 0.5340 & 1.7424 & 0.5365 & 1.8396 & 0.5315 & 1.6695 & 0.5294 & 1.7055 & 0.5203 & 1.6344 \\ 
       & & MSE & 0.0118 & 0.3706 & 0.0016 & 0.0993 & 0.0015 & 0.0400 & 0.0050 & 0.0175 & 0.0048 & 0.0105 \\
1.15 & 2 & Avg & 0.5175 & 1.5705 & 0.5196 & 1.6082 & 0.5154 & 1.5353 & 0.5136 & 1.5487 & 0.5057 & 1.5052 \\ 
       & & MSE & 0.0086 & 0.1633 & 0.0011 & 0.0269 & 0.0011 & 0.0200 & 0.0040 & 0.0082 & 0.0039 & 0.0073 \\ 
     & 3 & Avg & 0.5311 & 1.6775 & 0.5333 & 1.7381 & 0.5290 & 1.6267 & 0.5272 & 1.6504 & 0.5193 & 1.5973 \\ 
       & & MSE & 0.0100 & 0.2743 & 0.0013 & 0.0613 & 0.0012 & 0.0318 & 0.0043 & 0.0125 & 0.0041 & 0.0087 \\  
   \hline
\end{tabular}
\end{sidewaystable}
\begin{sidewaystable}
  \centering
\caption{Average values and the corresponding MSEs of the Bayes estimators with MCMC Method when $n=50$ and $m=40$.}
\label{MCMC:n50m40}
\begin{tabular}{ccccccccccccc}
\multicolumn{13}{c}{Informative Priors}\\
  \hline
 &&& \multicolumn{2}{c}{SEL} & \multicolumn{2}{c}{LINEX} & \multicolumn{2}{c}{LINEX} & \multicolumn{2}{c}{GEL} &\multicolumn{2}{c}{GEL}\\
 &&& & & \multicolumn{2}{c}{$\nu=-0.5$} & \multicolumn{2}{c}{$\nu=0.5$} & \multicolumn{2}{c}{$\kappa=-0.5$} &\multicolumn{2}{c}{$\kappa=0.5$}\\
 \cmidrule{4-13} 
 $T$& $R$&&$\aa$&$\bb$&$\aa$&$\bb$&$\aa$&$\bb$&$\aa$&$\bb$&$\aa$&$\bb$ \\ 
 \hline
     & 1 & Avg & 0.5178 & 1.5603 & 0.5186 & 1.5747 & 0.5170 & 1.5462 & 0.5163 & 1.5512 & 0.5132 & 1.5332 \\   
       & & MSE & 0.0017 & 0.0251 & 0.0002 & 0.0036 & 0.0002 & 0.0030 & 0.0008 & 0.0014 & 0.0008 & 0.0012 \\ 
0.21 & 2 & Avg & 0.5188 & 1.5599 & 0.5196 & 1.5755 & 0.5179 & 1.5448 & 0.5171 & 1.5502 & 0.5139 & 1.5306 \\  
       & & MSE & 0.0017 & 0.0241 & 0.0002 & 0.0035 & 0.0002 & 0.0028 & 0.0008 & 0.0013 & 0.0008 & 0.0011 \\  
     & 3 & Avg & 0.5186 & 1.5665 & 0.5195 & 1.5822 & 0.5178 & 1.5513 & 0.5171 & 1.5567 & 0.5139 & 1.5371 \\ 
       & & MSE & 0.0017 & 0.0252 & 0.0002 & 0.0037 & 0.0002 & 0.0030 & 0.0008 & 0.0014 & 0.0008 & 0.0012 \\
       \hline 
     & 1 & Avg & 0.5155 & 1.5481 & 0.5162 & 1.5581 & 0.5149 & 1.5382 & 0.5142 & 1.5417 & 0.5116 & 1.5290 \\ 
       & & MSE & 0.0017 & 0.0282 & 0.0002 & 0.0039 & 0.0002 & 0.0034 & 0.0008 & 0.0015 & 0.0008 & 0.0014 \\  
1.15 & 2 & Avg & 0.5166 & 1.5473 & 0.5173 & 1.5583 & 0.5158 & 1.5364 & 0.5152 & 1.5402 & 0.5124 & 1.5262 \\ 
       & & MSE & 0.0018 & 0.0264 & 0.0002 & 0.0036 & 0.0002 & 0.0032 & 0.0008 & 0.0014 & 0.0008 & 0.0013 \\ 
     & 3 & Avg & 0.5170 & 1.5497 & 0.5177 & 1.5609 & 0.5163 & 1.5388 & 0.5156 & 1.5427 & 0.5129 & 1.5285 \\ 
       & & MSE & 0.0019 & 0.0270 & 0.0002 & 0.0037 & 0.0002 & 0.0033 & 0.0009 & 0.0015 & 0.0008 & 0.0013 \\ 
   \hline
\end{tabular}
\begin{tabular}{ccccccccccccc}
\multicolumn{13}{c}{Non-Informative Priors}\\
  \hline
   &&& \multicolumn{2}{c}{SEL} & \multicolumn{2}{c}{LINEX} & \multicolumn{2}{c}{LINEX} & \multicolumn{2}{c}{GEL} &\multicolumn{2}{c}{GEL}\\
 &&& & & \multicolumn{2}{c}{$\nu=-0.5$} & \multicolumn{2}{c}{$\nu=0.5$} & \multicolumn{2}{c}{$\kappa=-0.5$} &\multicolumn{2}{c}{$\kappa=0.5$}\\
 \cmidrule{4-13} 
  $T$& $R$&&$\aa$&$\bb$&$\aa$&$\bb$&$\aa$&$\bb$&$\aa$&$\bb$&$\aa$&$\bb$ \\ 
 \hline
     & 1 & Avg & 0.5171 & 1.6144 & 0.5195 & 1.6768 & 0.5147 & 1.5622 & 0.5126 & 1.5839 & 0.5034 & 1.5247 \\  
       & & MSE & 0.0103 & 0.2243 & 0.0014 & 0.0425 & 0.0013 & 0.0260 & 0.0047 & 0.0115 & 0.0045 & 0.0091 \\
0.21 & 2 & Avg & 0.5186 & 1.6292 & 0.5215 & 1.7119 & 0.5158 & 1.5634 & 0.5134 & 1.5911 & 0.5028 & 1.5176 \\ 
       & & MSE & 0.0113 & 0.2788 & 0.0015 & 0.0597 & 0.0014 & 0.0317 & 0.0051 & 0.0144 & 0.0050 & 0.0110 \\  
     & 3 & Avg & 0.5202 & 1.6528 & 0.5229 & 1.7388 & 0.5175 & 1.5848 & 0.5151 & 1.6141 & 0.5050 & 1.5395 \\ 
       & & MSE & 0.0116 & 0.3012 & 0.0015 & 0.0658 & 0.0014 & 0.0341 & 0.0052 & 0.0154 & 0.0050 & 0.0114 \\
       \hline 
     & 1 & Avg & 0.5141 & 1.5574 & 0.5154 & 1.5768 & 0.5128 & 1.5388 & 0.5116 & 1.5459 & 0.5066 & 1.5230 \\ 
       & & MSE & 0.0058 & 0.0981 & 0.0007 & 0.0150 & 0.0007 & 0.0122 & 0.0026 & 0.0048 & 0.0026 & 0.0044 \\
1.15 & 2 & Avg & 0.5134 & 1.5415 & 0.5149 & 1.5630 & 0.5119 & 1.5209 & 0.5106 & 1.5284 & 0.5049 & 1.5020 \\ 
       & & MSE & 0.0060 & 0.0874 & 0.0008 & 0.0126 & 0.0007 & 0.0107 & 0.0028 & 0.0046 & 0.0028 & 0.0043 \\ 
     & 3 & Avg & 0.5144 & 1.5473 & 0.5158 & 1.5692 & 0.5130 & 1.5264 & 0.5117 & 1.5341 & 0.5063 & 1.5074 \\ 
       & & MSE & 0.0061 & 0.0915 & 0.0008 & 0.0134 & 0.0008 & 0.0113 & 0.0028 & 0.0048 & 0.0028 & 0.0044 \\  
   \hline
\end{tabular}
\end{sidewaystable}
\begin{table}
  \centering
\caption{Average values and the corresponding MSEs of the SPT estimators when $n=30$.}
\label{SPT:n30}
\begin{tabular}{ccccccccccc}
\multicolumn{11}{c}{$m=15$}\\
  \hline
 &&& \multicolumn{2}{c}{NR} & \multicolumn{2}{c}{SEL} & \multicolumn{2}{c}{LINEX} & \multicolumn{2}{c}{GEL} \\
 \cmidrule{4-11} 
 $T$& $R$&&$\aa$&$\bb$&$\aa$&$\bb$&$\aa$&$\bb$&$\aa$&$\bb$ \\ 
 \hline
     & 1 & Avg & 0.5803 & 1.7754 & 0.5894 & 1.6347 & 0.5886 & 1.6257 & 0.5842 & 1.6175 \\   
       & & MSE & 0.0252 & 0.6967 & 0.0114 & 0.0237 & 0.0112 & 0.0212 & 0.0107 & 0.0195 \\ 
0.21 & 2 & Avg & 0.5694 & 1.7186 & 0.5935 & 1.6357 & 0.5926 & 1.6256 & 0.5885 & 1.6163 \\  
       & & MSE & 0.0261 & 0.3938 & 0.0117 & 0.0224 & 0.0116 & 0.0198 & 0.0110 & 0.0177 \\  
     & 3 & Avg & 0.5693 & 1.6991 & 0.5878 & 1.6349 & 0.5869 & 1.6254 & 0.5826 & 1.6167 \\ 
       & & MSE & 0.0242 & 0.3782 & 0.0113 & 0.0225 & 0.0112 & 0.0200 & 0.0107 & 0.0181 \\
       \hline 
     & 1 & Avg & 0.5884 & 1.8161 & 0.5872 & 1.6363 & 0.5863 & 1.6275 & 0.5818 & 1.6196 \\ 
       & & MSE & 0.0253 & 0.5936 & 0.0112 & 0.0251 & 0.0111 & 0.0227 & 0.0105 & 0.0209 \\  
1.15 & 2 & Avg & 0.5532 & 1.6404 & 0.5735 & 1.6302 & 0.5725 & 1.6214 & 0.5671 & 1.6130 \\ 
       & & MSE & 0.0187 & 0.1543 & 0.0100 & 0.0236 & 0.0099 & 0.0214 & 0.0094 & 0.0198 \\ 
     & 3 & Avg & 0.5752 & 1.7205 & 0.5760 & 1.6317 & 0.5750 & 1.6232 & 0.5702 & 1.6151 \\ 
       & & MSE & 0.0209 & 0.3122 & 0.0103 & 0.0251 & 0.0102 & 0.0228 & 0.0097 & 0.0214 \\ 
   \hline
\end{tabular}
\begin{tabular}{ccccccccccc}
\multicolumn{11}{c}{$m=20$}\\
  \hline
 &&& \multicolumn{2}{c}{NR} & \multicolumn{2}{c}{SEL} & \multicolumn{2}{c}{LINEX} & \multicolumn{2}{c}{GEL} \\
 \cmidrule{4-11} 
 $T$& $R$&&$\aa$&$\bb$&$\aa$&$\bb$&$\aa$&$\bb$&$\aa$&$\bb$ \\ 
 \hline
     & 1 & Avg & 0.5677 & 1.6594 & 0.5842 & 1.6313 & 0.5833 & 1.6227 & 0.5788 & 1.6148 \\  
       & & MSE & 0.0223 & 0.2083 & 0.0111 & 0.0226 & 0.0110 & 0.0204 & 0.0104 & 0.0187 \\
0.21 & 2 & Avg & 0.5664 & 1.6919 & 0.5879 & 1.6368 & 0.5869 & 1.6272 & 0.5827 & 1.6185 \\ 
       & & MSE & 0.0243 & 0.2876 & 0.0113 & 0.0233 & 0.0112 & 0.0207 & 0.0107 & 0.0188 \\  
     & 3 & Avg & 0.5670 & 1.6887 & 0.5875 & 1.6347 & 0.5866 & 1.6252 & 0.5819 & 1.6162 \\ 
       & & MSE & 0.0234 & 0.2773 & 0.0113 & 0.0228 & 0.0111 & 0.0203 & 0.0106 & 0.0185 \\
       \hline 
     & 1 & Avg & 0.5704 & 1.6809 & 0.5721 & 1.6296 & 0.5710 & 1.6218 & 0.5658 & 1.6142 \\ 
       & & MSE & 0.0186 & 0.1565 & 0.0098 & 0.0257 & 0.0097 & 0.0237 & 0.0092 & 0.0225 \\
1.15 & 2 & Avg & 0.5549 & 1.6182 & 0.5687 & 1.6260 & 0.5677 & 1.6176 & 0.5622 & 1.6098 \\ 
       & & MSE & 0.0168 & 0.1237 & 0.0096 & 0.0241 & 0.0095 & 0.0222 & 0.0089 & 0.0210 \\ 
     & 3 & Avg & 0.5490 & 1.6052 & 0.5628 & 1.6222 & 0.5616 & 1.6140 & 0.5563 & 1.6066 \\ 
       & & MSE & 0.0159 & 0.1309 & 0.0090 & 0.0234 & 0.0089 & 0.0216 & 0.0084 & 0.0202 \\  
   \hline
\end{tabular}
\end{table}
\begin{table}
  \centering
\caption{Average values and the corresponding MSEs of the SPT estimators when $n=50$.}
\label{SPT:n50}
\begin{tabular}{ccccccccccc}
\multicolumn{11}{c}{$m=25$}\\
  \hline
 &&& \multicolumn{2}{c}{NR} & \multicolumn{2}{c}{SEL} & \multicolumn{2}{c}{LINEX} & \multicolumn{2}{c}{GEL} \\
 \cmidrule{4-11} 
 $T$& $R$&&$\aa$&$\bb$&$\aa$&$\bb$&$\aa$&$\bb$&$\aa$&$\bb$ \\ 
 \hline
     & 1 & Avg & 0.5554 & 1.6663 & 0.5635 & 1.6327 & 0.5626 & 1.6250 & 0.5578 & 1.6178 \\   
       & & MSE & 0.0167 & 0.1846 & 0.0089 & 0.0244 & 0.0088 & 0.0223 & 0.0083 & 0.0210 \\ 
0.21 & 2 & Avg & 0.5530 & 1.6441 & 0.5769 & 1.6334 & 0.5761 & 1.6242 & 0.5714 & 1.6159 \\  
       & & MSE & 0.0189 & 0.2094 & 0.0103 & 0.0228 & 0.0102 & 0.0204 & 0.0097 & 0.0187 \\  
     & 3 & Avg & 0.5491 & 1.6271 & 0.5643 & 1.6308 & 0.5633 & 1.6227 & 0.5587 & 1.6150 \\ 
       & & MSE & 0.0158 & 0.1365 & 0.0090 & 0.0233 & 0.0089 & 0.0213 & 0.0084 & 0.0200 \\
       \hline 
     & 1 & Avg & 0.5617 & 1.6945 & 0.5622 & 1.6347 & 0.5610 & 1.6268 & 0.5562 & 1.6201 \\ 
       & & MSE & 0.0164 & 0.1755 & 0.0087 & 0.0258 & 0.0086 & 0.0238 & 0.0081 & 0.0224 \\  
1.15 & 2 & Avg & 0.5394 & 1.6026 & 0.5505 & 1.6234 & 0.5494 & 1.6158 & 0.5445 & 1.6084 \\ 
       & & MSE & 0.0128 & 0.0895 & 0.0075 & 0.0242 & 0.0074 & 0.0225 & 0.0070 & 0.0214 \\ 
     & 3 & Avg & 0.5557 & 1.6564 & 0.5554 & 1.6270 & 0.5544 & 1.6196 & 0.5496 & 1.6127 \\ 
       & & MSE & 0.0144 & 0.1185 & 0.0079 & 0.0261 & 0.0078 & 0.0243 & 0.0074 & 0.0231 \\ 
   \hline
\end{tabular}
\begin{tabular}{ccccccccccc}
\multicolumn{11}{c}{$m=40$}\\
  \hline
 &&& \multicolumn{2}{c}{NR} & \multicolumn{2}{c}{SEL} & \multicolumn{2}{c}{LINEX} & \multicolumn{2}{c}{GEL} \\
 \cmidrule{4-11} 
 $T$& $R$&&$\aa$&$\bb$&$\aa$&$\bb$&$\aa$&$\bb$&$\aa$&$\bb$ \\ 
 \hline
     & 1 & Avg & 0.5414 & 1.6053 & 0.5566 & 1.6278 & 0.5557 & 1.6203 & 0.5512 & 1.6136 \\  
       & & MSE & 0.0142 & 0.1050 & 0.0081 & 0.0227 & 0.0080 & 0.0209 & 0.0076 & 0.0196 \\
0.21 & 2 & Avg & 0.5463 & 1.6119 & 0.5659 & 1.6287 & 0.5649 & 1.6207 & 0.5598 & 1.6132 \\ 
       & & MSE & 0.0156 & 0.1253 & 0.0091 & 0.0225 & 0.0090 & 0.0205 & 0.0085 & 0.0191 \\  
     & 3 & Avg & 0.5477 & 1.6250 & 0.5632 & 1.6323 & 0.5623 & 1.6244 & 0.5576 & 1.6170 \\ 
       & & MSE & 0.0158 & 0.1297 & 0.0089 & 0.0232 & 0.0088 & 0.0211 & 0.0083 & 0.0197 \\
       \hline 
     & 1 & Avg & 0.5277 & 1.5916 & 0.5288 & 1.6116 & 0.5278 & 1.6053 & 0.5237 & 1.5994 \\ 
       & & MSE & 0.0085 & 0.0606 & 0.0044 & 0.0253 & 0.0043 & 0.0243 & 0.0040 & 0.0238 \\
1.15 & 2 & Avg & 0.5305 & 1.5866 & 0.5352 & 1.6138 & 0.5344 & 1.6075 & 0.5298 & 1.6012 \\ 
       & & MSE & 0.0095 & 0.0563 & 0.0054 & 0.0243 & 0.0053 & 0.0231 & 0.0049 & 0.0224 \\ 
     & 3 & Avg & 0.5304 & 1.5906 & 0.5340 & 1.6166 & 0.5332 & 1.6100 & 0.5290 & 1.6039 \\ 
       & & MSE & 0.0094 & 0.0592 & 0.0052 & 0.0242 & 0.0051 & 0.0231 & 0.0048 & 0.0223  \\  
   \hline
\end{tabular}
\end{table}
\begin{sidewaystable}
\centering
\caption{Confidence intervals and covarage probabilities of NR and MCMC methods, when $n=30$. (U:upper, L:lower, IL: interval length, CP: coverage probability}
\begin{tabular}{ccccccccccccccccc} \label{CI:CP:n30}
 $T$& $R$ && L & U & IL & CP && L & U & IL & CP && L & U & IL & CP \\ 
  \hline
  \multicolumn{17}{c}{$m=15$}\\
  & &&\multicolumn{4}{c}{NR} && \multicolumn{4}{c}{MCMC:Informative} && \multicolumn{4}{c}{MCMC:Non-Informative}\\
  \cmidrule{4-7}\cmidrule{9-12} \cmidrule{14-17} 
      & 1 & $\aa$ & 0.2756 & 0.8127 & 0.5371 & 95.48 && 0.4054 & 0.6486 & 0.2432 & 99.86 & & 0.2990 & 0.8151 & 0.5161 & 93.60 \\ 
        & & $\bb$ & 0.3339 & 3.2387 & 2.9048 & 94.26 && 1.1011 & 2.1334 & 1.0323 & 99.80 & & 0.7430 & 3.5039 & 2.7609 & 92.42 \\ 
 0.21 & 2 & $\aa$ & 0.2289 & 0.8356 & 0.6066 & 94.76 && 0.4009 & 0.6560 & 0.2551 & 99.96 & & 0.2624 & 0.8596 & 0.5972 & 95.40 \\ 
        & & $\bb$ & 0.0442 & 3.4998 & 3.4555 & 94.04 && 1.0715 & 2.1728 & 1.1014 & 99.94 & & 0.5673 & 3.9731 & 3.4058 & 93.40 \\ 
      & 3 & $\aa$ & 0.2587 & 0.8072 & 0.5486 & 95.00 && 0.4041 & 0.6506 & 0.2464 & 99.96 & & 0.2825 & 0.8161 & 0.5336 & 94.74 \\
\hline
\end{tabular}
\end{sidewaystable}

\begin{sidewaystable}
\centering
\caption{Confidence intervals and covarage probabilities of NR and MCMC methods, when $n=50$. (U:upper, L:lower, IL: interval length, CP: coverage probability)}
\begin{tabular}{ccccccccccccccccc} \label{CI:CP:n50}
 $T$& $R$ && L & U & IL & CP && L & U & IL & CP && L & U & IL & CP \\ 
  \hline
  \multicolumn{17}{c}{$m=25$}\\
  & &&\multicolumn{4}{c}{NR} && \multicolumn{4}{c}{MCMC:Informative} && \multicolumn{4}{c}{MCMC:Non-Informative}\\
  \cmidrule{4-7}\cmidrule{9-12} \cmidrule{14-17} 
      & 1 & $\aa$ & 0.3302 & 0.7297 & 0.3995 & 94.54 && 0.4135 & 0.6342 & 0.2207 & 99.72 && 0.3468 & 0.7409 & 0.3941 & 93.44 \\ 
        & & $\bb$ & 0.6798 & 2.7055 & 2.0257 & 94.78 && 1.1345 & 2.0849 & 0.9504 & 99.90 && 0.9079 & 2.9515 & 2.0436 & 92.80 \\ 
 0.21 & 2 & $\aa$ & 0.2860 & 0.7589 & 0.4729 & 95.02 && 0.4074 & 0.6452 & 0.2378 & 99.82 && 0.3143 & 0.7852 & 0.4709 & 94.70 \\ 
        & & $\bb$ & 0.4227 & 2.9116 & 2.4889 & 94.18 && 1.0936 & 2.1339 & 1.0403 & 99.86 && 0.7368 & 3.3120 & 2.5752 & 94.20 \\ 
      & 3 & $\aa$ & 0.3214 & 0.7253 & 0.4039 & 95.52 && 0.4125 & 0.6368 & 0.2243 & 99.84 && 0.3394 & 0.7411 & 0.4017 & 94.72 \\
\hline
\end{tabular}
\end{sidewaystable}
\begin{table}[ht]
\centering
\caption{Estimation values of listed methods for Carbon Fibre data}\label{app:est}
\begin{tabular}{lcccccccc}
  & \multicolumn{2}{c}{Sch 1}&&\multicolumn{2}{c}{Sch 2}&&\multicolumn{2}{c}{Sch 3}\\ \cmidrule{2-3}\cmidrule{5-6}\cmidrule{8-9}
MLE Method & $\aa$ & $\bb$ &&$\aa$&$\bb$&&$\aa$&$\bb$\\ 
  \hline
  NR  & 2.2542 & 0.3980 && 2.2276 & 0.3441 && 2.1092 & 0.3896\\ 
  EM  & 2.4056 & 0.3577 && 2.3619 & 0.3049 && 2.1217 & 0.3860\\ 
  SEM & 1.7966 & 0.4142 && 1.5622 & 0.3637 && 2.3277 & 0.4440\\ 
  \hline
  Tierney-Kadane Method \\
  SEL               & 2.2505 & 0.3991 && 2.2276 & 0.3467 && 2.1029 & 0.3911\\ 
  LINEX($\nu=-0.5$) & 2.2764 & 0.4005 && 2.2559 & 0.3485 && 2.1227 & 0.3925\\ 
  LINEX($\nu=0.5$)  & 2.2260 & 0.3978 && 2.2018 & 0.3451 && 2.0841 & 0.3897\\ 
  GEL($\kappa=-0.5$)& 2.2393 & 0.3958 && 2.2155 & 0.3418 && 2.0937 & 0.3875\\ 
  GEL($\kappa=0.5$) & 2.2169 & 0.3890 && 2.1912 & 0.3321 && 2.0752 & 0.3803\\ 
  \hline 
  MCMC Method \\
  SEL               & 2.2579 & 0.3975 && 2.3711 & 0.3341 && 2.0818 & 0.3977\\ 
  LINEX($\nu=-0.5$) & 2.2898 & 0.3990 && 2.4008 & 0.3357 && 2.0980 & 0.3991\\ 
  LINEX($\nu=0.5$)  & 2.2265 & 0.3960 && 2.3420 & 0.3325 && 2.0653 & 0.3963\\ 
  GEL($\kappa=-0.5$)& 2.2438 & 0.3938 && 2.3587 & 0.3295 && 2.0737 & 0.3943\\ 
  GEL($\kappa=0.5$) & 2.2150 & 0.3865 && 2.3335 & 0.3203 && 2.0570 & 0.3873\\ 
  \hline
  Shrinkage Method \\
  NR               & 2.2524 & 0.3985 && 2.2276 & 0.3454 && 2.1060 & 0.3904\\ 
  SEL              & 2.2505 & 0.3991 && 2.2276 & 0.3467 && 2.1029 & 0.3911\\ 
  LINEX($\nu=0.5$) & 2.2634 & 0.3998 && 2.2418 & 0.3476 && 2.1128 & 0.3918\\ 
  GEL($\kappa=0.5$)& 2.2449 & 0.3974 && 2.2216 & 0.3442 && 2.0983 & 0.3893\\ 
   \hline
\end{tabular}
\end{table}

\begin{table}[ht]
\centering
\caption{Confident intervals and interval lengths of NR and MCMC methods for Carbon Fibre data (U:upper, L:lower, IL: interval length)} \label{app:CI}
\begin{tabular}{clccccccc}
  && \multicolumn{3}{c}{$\alpha$}&&\multicolumn{3}{c}{$\bb$}  \\ \cmidrule{3-5}\cmidrule{7-9}
Scheme & Method & L & U & IL & &L & U & IL \\ 
  \hline
1&NR   & 1.6321 & 2.8764 & 1.2443 && 0.2539 & 0.5420 & 0.2880 \\ 
&MCMC & 1.5452 & 3.0372 & 1.4919 && 0.2658 & 0.5615 & 0.2957 \\  
   \hline
2&NR   & 1.5815 & 2.8737 & 1.2923 && 0.1825 & 0.5056 & 0.3231 \\ 
 &MCMC & 1.7059 & 3.0818 & 1.3759 && 0.1989 & 0.5192 & 0.3203 \\ 
   \hline
3&NR   & 1.5641 & 2.6542 & 1.0902 && 0.2426 & 0.5367 & 0.2941 \\ 
 &MCMC & 1.5316 & 2.5642 & 1.0325 && 0.2597 & 0.5609 & 0.3012 \\ 
   \hline
\end{tabular}
\end{table}

\begin{thebibliography}{99}
\bibitem[Ahmed and Saleh, 1990]{ahmed1990}
Ahmed, S. E., Saleh, A. M. E. (1990). Estimation strategies for the intercept vector in a simple linear multivariate normal regression model. Computational Statistics \& Data Analysis, 10(3), 193--206.

\bibitem[Asgharzadeh et al., 2015]{Asgar}
Asgharzadeh, A., Valiollahi, R., Kundu, D. (2015). Prediction for future failures in Weibull distribution under hybrid censoring. Journal of Statistical Computation and Simulation, 85(4), 824--838.

\bibitem[Asl et al., 2018]{Reza2018}
Asl, M. N., Belaghi, R. A., Bevrani, H. (2018). Classical and Bayesian inferential approaches using Lomax model under progressively type-I hybrid censoring. Journal of Computational and Applied Mathematics, https://doi.org/10.1016/j.cam.2018.04.028.


\bibitem[Bader and Priest, 1982]{real_data}
Bader, M. G., Priest, A. M. (1982). Statistical aspects of fibre and bundle strength in hybrid composites. Progress in science and engineering of composites, 1129--1136.

\bibitem[Balakrishnan, 2007]{Bala2007prog}
Balakrishnan, N. (2007). Progressive censoring methodology: an appraisal. Test, 16(2), 211.

\bibitem[Balakrishnan and
Aggarwala, 2000]{Bala2000}
Balakrishnan, N., Aggarwala, R. (2000). Progressive censoring: theory, methods, and applications, Springer, New York.

\bibitem[Balakrishnan and Kateri, 2008]{Kateri2008}
Balakrishnan, N., Kateri, M. (2008). On the maximum likelihood estimation of parameters of Weibull distribution based on complete and censored data. Statistics \& Probability Letters, 78(17), 2971--2975.

\bibitem[Balakrishnan and Kundu, 2013]{BK2013}
Balakrishnan, N., Kundu, D. (2013). Hybrid censoring: models, inferential results and applications. Computational Statistics \& Data Analysis, 57(1), 166--209.

\bibitem[Banerjee and Kundu, 2008]{BK2008}
Banerjee, A., Kundu, D. (2008). Inference based on type-II hybrid censored data from a Weibull distribution. IEEE Transactions on reliability, 57(2), 369--378.

\bibitem[Bancroft, 1944]{bancroft}
Bancroft, T. A. (1944). On biases in estimation due to the use of preliminary tests of significance. The Annals of Mathematical Statistics, 15(2), 190--204.


\bibitem[Chen and Bhattacharya, 1988]{ChenBhat1988}
Chen, S., Bhattacharya, G.K., 1988. Exact confidence bounds for an exponential parameter under hybrid censoring. Communications in Statistics--Theory
and Methods 17, 1857--1870.

\bibitem[Chen and Shao, 1999]{chen-shao}
Chen, M. H., Shao, Q. M. (1999). Monte Carlo estimation of Bayesian credible and HPD intervals. Journal of Computational and Graphical Statistics, 8(1), 69--92.

\bibitem[Childs et al., 2008]{childs2008}
Childs A., Chandrasekar B., Balakrishnan N. (2008) Exact Likelihood Inference for an Exponential Parameter Under Progressive Hybrid Censoring Schemes. In: Vonta F., Nikulin M., Limnios N., Huber-Carol C. (eds) Statistical Models and Methods for Biomedical and Technical Systems. Statistics for Industry and Technology. Birkhäuser Boston.

\bibitem[Diebolt and Celeux, 1993]{EM1993}
Diebolt, J., Celeux, G. (1993). Asymptotic properties of a stochastic EM algorithm for estimating mixing proportions. Stochastic Models, 9(4), 599-613.

\bibitem[Dempster et al., 1977]{em1977}
Dempster, A. P., Laird, N. M., Rubin, D. B. (1977). Maximum likelihood from incomplete data via the EM algorithm. Journal of the royal statistical society. Series B (methodological), 1--38.


\bibitem[Dey et al., 2016]{Dey2016b}
Dey, S., Singh, S., Tripathi, Y. M., Asgharzadeh, A. (2016). Estimation and prediction for a progressively censored generalized inverted exponential distribution. Statistical Methodology, 32, 185--202.

\bibitem[Epstein, 1954]{epstein1954}
Epstein, B. (1954). Truncated life tests in the exponential case, The Annals of Mathematical Statistics,
25(3), 555--564.

\bibitem[Jia et al., 2018]{jia}
Jia, X., Nadarajah, S., Guo, B. (2018). Exact Inference on Weibull Parameters With Multiply Type-I Censored Data. IEEE Transactions on Reliability, 67(2), 432--445.


\bibitem[Huang and Wu, 2012]{HuangWu2012}
Huang, S. R., Wu, S. J. (2012). Bayesian estimation and prediction for Weibull model with progressive censoring. Journal of Statistical Computation and Simulation, 82(11), 1607--1620.

\bibitem[Judge and Bock, 1978]{JB}
Judge, G. G., Bock, M. E. (1978) The statistical implicatinos of pre-test and stein-rule estimators in econometrics.
North-Holland, Amsterdam

\bibitem[Kibria and Saleh, 2010]{kibria2010}
Kibria, B. G., Saleh, A. M. E. (2010). Preliminary test estimation of the parameters of exponential and Pareto distributions for censored samples. Statistical Papers, 51(4), 757--773.

\bibitem[Kim et al., 2011]{Kim2011}
Kim, C., Jung, J., Chung, Y. (2011). Bayesian estimation for the exponentiated Weibull model under Type II progressive censoring. Statistical Papers, 52(1), 53--70.

\bibitem[Kundu and Joarder, 2006]{kundu2006}
Kundu, D., Joarder, A. (2006).  Analysis of type-II progressively hybrid censored data, Computational
Statistics Data Analysis,  50(10), 2509--2528.

\bibitem[Kundu, 2008]{Kundu2008}
Kundu, D. (2008). Bayesian inference and life testing plan for the Weibull distribution in presence of progressive censoring. Technometrics, 50(2), 144--154.

\bibitem[Lawless, 2003]{Lawless}
Lawless, J .F. (2003). Statistical models and methods for lifetime data, 2nd~ed. Hoboken: Wiley.

\bibitem[Lin et al., 2009]{lin2009}
Lin, C. T., Ng, H. K. T., Chan, P. S. (2009). Statistical inference of Type-II progressively hybrid censored data with Weibull lifetimes. Communications in Statistics--Theory and Methods, 38(10), 1710--1729.

\bibitem[Lin and Huang, 2012]{LinHuang2012}  
Lin, C. T., Huang, Y. L. (2012). On progressive hybrid censored exponential distribution. Journal of Statistical Computation and Simulation, 82(5), 689--709.

\bibitem[Lin et al., 2012]{Lin2012}
Lin, C. T., Chou, C. C., Huang, Y. L. (2012). Inference for the Weibull distribution with progressive hybrid censoring. Computational Statistics \& Data Analysis, 56(3), 451--467.



\bibitem[Louis, 1982]{Louis1982}
Louis, T. A. (1982). Finding the observed information matrix when using the em algorithm, Journal of the Royal Statistical Society. Series B (Methodological), 226--233.
 
 \bibitem[Mokhtari et al., 2011]{moktari}
 Mokhtari, E. B., Rad, A. H., Yousefzadeh, F. (2011). Inference for Weibull distribution based on progressively Type-II hybrid censored data. Journal of Statistical Planning and Inference, 141(8), 2824--2838.
 


\bibitem[Ng et al., 2002]{Ng2002}
Ng, H. K. T., Chan, P. S., Balakrishnan, N. (2002). Estimation of parameters from progressively censored data using EM algorithm. Computational Statistics \& Data Analysis, 39(4), 371--386.

\bibitem[Pandey, 1983]{Pandey1983}
Pandey, M. (1983). Shrunken estimators of Weibull shape parameter in censored samples. IEEE Transactions on Reliability, 32(2), 200--203.

\bibitem[Pandey and Singh, 1993]{PS1993}
Pandey M., Singh U. S. (1993). Shrunken estimators of Weibull shape parameter from Type-II censored
samples. IEEE Transactions on Reliability 42, 81--86

\bibitem[R Core Team, 2018]{R2018}
R Core Team (2018). R: A language and environment for statistical computing. R Foundation for Statistical Computing, Vienna, Austria. URL https://www.R-project.org/.

   
\bibitem[Saleh, 2006]{Saleh2006}
Saleh, A. M. E. (2006) Theory of preliminary test and stein-type estimations with applications. Wiley,
New York.

\bibitem[Saleh and Kibria, 1993]{SalehKibria1993}
Saleh, A. M. E., Kibria, B. M. G. (1993). Performance of some new preliminary test ridge regression estimators and their properties. Communications in Statistics--Theory and Methods, 22(10), 2747--2764.

\bibitem[Singh and Bhatkulikar, 1978]{SB1978}
Singh, J., Bhatkulikar, S. G. (1978). Shrunken estimation in Weibull distribution. Sankhyā: The Indian Journal of Statistics, Series B, 382--393.

\bibitem[Singh and Shukla, 2000]{SS2000}
Singh, H. P., Shukla, S. K. (2000). Estimation in the Two-parameter Weibull distribution with Prior Information. IAPQR TRANSACTIONS, 25(2), 107--118.


\bibitem[Thompson, 1968]{tompson}
Thompson, J. R. (1968). Some shrinkage techniques for estimating the mean. Journal of the American Statistical Association, 63(321), 113--122.

\bibitem[Tierney and Kadane, 1986]{Tierney}
Tierney, L., Kadane, J. B. (1986).  Accurate approximations for posterior moments and marginal densities, Journal of the American Statistical Association, 81(393), 82--86.

\bibitem[Varian, 1975]{linex}
Varian, H. R. (1975). A Bayesian approach to real estate assessment. Studies in Bayesian Econometric and Statistics in honor of Leonard J. Savage, 195--208.

\bibitem[Wei and Tanner, 1990]{WeiTanner1990}
Wei, G. C., Tanner, M. A. (1990). A Monte Carlo implementation of the EM algorithm and the poor man's data augmentation algorithms. Journal of the American Statistical Association, 85(411), 699--704.

\end{thebibliography}
\end{document}